\documentclass[11pt,leqno]{article}
\usepackage[mathscr]{eucal}
\usepackage{amsmath,amssymb,latexsym,theorem,bbm}
\usepackage{color}
\usepackage{appendix}

\newcommand{\NN}{\mathbb{N}}

\newcommand{\RR}{\mathbb{R}}

\newcommand{\ZZ}{\mathbb{Z}}

\newcommand{\cA}{{\mathcal A}}
\newcommand{\cB}{{\mathcal B}}

\newcommand{\cC}{{\mathcal C}}

\newcommand{\cF}{{\mathcal F}}

\newcommand{\cL}{{\mathcal L}}

\newcommand{\cN}{{\mathcal N}}

\newcommand{\dd}{\mathrm{d}}
\newcommand{\ee}{\mathrm{e}}

\DeclareMathOperator*{\argmin}{arg\,min}

\newcommand{\EE}{\operatorname{\mathbb{E}}}
\newcommand{\PP}{\operatorname{\mathbb{P}}}

\newcommand{\diag}{\operatorname{diag}}

\renewcommand{\mid}{\,|\,}

\renewcommand{\leq}{\leqslant}
\renewcommand{\geq}{\geqslant}

\newcommand{\stoch}{\stackrel{\PP}{\longrightarrow}}
\newcommand{\distr}{\stackrel{\cL}{\longrightarrow}}

\newcommand{\proofend}{\hfill\mbox{$\Box$}}

\numberwithin{equation}{section}

\theoremstyle{change} \theorembodyfont{\em}
\newtheorem{Lem}{Lemma.}[section]
\newtheorem{Thm}[Lem]{Theorem.}
\newtheorem{Pro}[Lem]{Proposition.}

\newtheorem{Def}[Lem]{Definition.}

\theorembodyfont{\rm}
\newtheorem{Rem}[Lem]{Remark.}

\begin{document}

\begin{center}
 {\bfseries\Large
  Parameter estimation for a subcritical affine two factor model} \\[5mm]
 {\sc\large
  M\'aty\'as $\text{Barczy}^{*}$,
  \ Leif $\text{D\"oring}$,
  \ Zenghu $\text{Li}$,
  \ Gyula $\text{Pap}$}
\end{center}


\renewcommand{\thefootnote}{}
\footnote{$*$ Corresponding author}
\footnote{\textit{2010 Mathematics Subject Classifications\/}:
          62F12, 60J25.}
\footnote{\textit{Key words and phrases\/}:
 affine process, maximum likelihood estimator, least squares estimator.}
\vspace*{0.2cm}
\footnote{
The research of M. Barczy and G. Pap was realized in the frames of T\'AMOP 4.2.4.\ A/2-11-1-2012-0001 ,,National Excellence
Program -- Elaborating and operating an inland student and researcher personal support
system''. The project was subsidized by the European Union and co-financed by the
European Social Fund.
Z. Li has been partially supported by NSFC under Grant No.\ 11131003 and 973
 Program under Grant No.\ 2011CB808001.}

\vspace*{-10mm}

\begin{abstract}
For an affine two factor model, we study the asymptotic properties of the maximum likelihood and
 least squares estimators of some appearing parameters in the so-called subcritical (ergodic) case
 based on continuous time observations.
We prove strong consistency and asymptotic normality of the estimators in question.
\end{abstract}

\section{Introduction}

We consider the following 2-dimensional affine process (affine two factor model)
 \begin{align}\label{2dim_affine}
  \begin{cases}
   \dd Y_t = (a-bY_t)\,\dd t + \sqrt{Y_t}\,\dd L_t,\\
   \dd X_t = (m-\theta X_t)\,\dd t + \sqrt{Y_t}\,\dd B_t,
  \end{cases}
  t\geq 0,
 \end{align}
 where \ $a>0$, \ $b, m, \theta \in \RR$, \ and \ $(L_t)_{t\geq 0}$ \ and
 \ $(B_t)_{t\geq 0}$ \ are independent standard Wiener processes.
Note that the process \ $(Y_t)_{t\geq 0}$ \ given by the first SDE of
 \eqref{2dim_affine} is the so-called Cox--Ingersol--Ross (CIR) process which
 is a continuous state branching process with branching mechanism
 \ $b z + z^2 / 2$, \ $z \geq 0$, \ and with immigration mechanism \ $a z$,
 $z \geq 0$.
\ Chen and Joslin \cite{CheJos} applied \eqref{2dim_affine} for modelling
 quantitative impact of stochastic recovery on the pricing of defaultable
 bonds, see their equations (25) and (26).

The process \ $(Y, X)$ \ given by \eqref{2dim_affine} is a special affine
 diffusion process.
The set of affine processes contains a large class of important Markov processes such as
 continuous state branching processes and Ornstein--Uhlenbeck processes.
Further, a lot of models in financial mathematics are affine
 such as the Heston model \cite{Hes}, the model of Barndorff-Nielsen
 and Shephard \cite{BarShe} or the model due to Carr and Wu \cite{CarWu}.
A precise mathematical formulation and a complete characterization of regular
 affine processes are due to Duffie et al. \cite{DufFilSch}.
These processes are widely used in financial mathematics due to their
 computational tractability, see Gatheral \cite{Gat}.

This article is devoted to estimate the parameters \ $a$, $b$, $m$ \ and \ $\theta$ \ from some continuously observed real data set
 \ $(Y_t,X_t)_{t\in[0,T]}$, \ where \ $T>0$.
\ To the best knowledge of the authors the parameter estimation problem for multi-dimensional affine processes has not been tackled so far.
Since affine processes are frequently used in financial mathematics, the question of parameter estimation for them
 is of high importance.
In Barczy et al.~\cite{BarDorLiPap} we started the discussion with a simple non-trivial two-dimensional affine diffusion process
 given by \eqref{2dim_affine} in the so called critical case: \ $b\geq 0$, $\theta=0$ \ or \ $b=0$, $\theta\geq 0$
 \ (for the definition of criticality, see Section \ref{Section_preliminaires}).
In the special critical case \ $b=0$, $\theta=0$ \ we described the asymptotic behavior of least squares
 estimator (LSE) of \ $(m,\theta)$ \ from some discretely observed low frequency real data set \ $X_0,X_1,\ldots,X_n$ \ as \ $n\to\infty$.
\ The description of the asymptotic behavior of the LSE of \ $(m,\theta)$ \ in the other critical cases
 \ $b=0$, $\theta>0$ \ or \ $b>0$, $\theta=0$ \ remained opened.
\ In this paper we deal with the same model \eqref{2dim_affine} but in the so-called subcritical (ergodic) case:
 \ $b>0$, $\theta>0$, \ and we consider the maximum likelihood estimator (MLE) of \ $(a,b,m,\theta)$ \
 using some continuously observed real data set \ $(Y_t,X_t)_{t\in[0,T]}$, \ where \ $T>0$, \
 and the LSE of \ $(m,\theta)$ \ using some continuously observed real data set \ $(X_t)_{t\in[0,T]}$, \ where \ $T>0$.
\ For studying the asymptotic behaviour of the MLE and LSE in the subcritical (ergodic) case,
 one first needs to examine the question of existence of a unique stationary distribution
 and ergodicity for the model given by \eqref{2dim_affine}.
In a companion paper Barczy et al.~\cite{BarDorLiPap2} we solved this problem,
 see also Theorem \ref{Thm_ergodic12}.
Further, in a more general setup by replacing the CIR process \ $(Y_t)_{t\geq 0}$ \ in the first SDE of \eqref{2dim_affine}
 by a so-called \ $\alpha$-root process (stable CIR process) with \ $\alpha\in(1,2)$, \ the existence of a unique
 stationary distribution for the corresponding model was proved in Barczy et al. \cite{BarDorLiPap2}.

In general, parameter estimation for subcritical (also called ergodic) models has a long history, see, e.g.,
 the monographs of Liptser and Shiryaev \cite[Chapter 17]{LipShiII},
 Kutoyants \cite{Kut}, Bishwal \cite{Bis} and the papers of Klimko and Nelson \cite{KliNel} and S{\o}rensen \cite{Sor}.
In case of the one-dimensional CIR process \ $Y$, \ the parameter estimation of \ $a$ \ and \ $b$ \ goes back to
 Overbeck and Ryd\'en \cite{OveRyd} (conditional LSE), Overbeck \cite{Ove} (MLE), and see also Bishwal \cite[Example 7.6]{Bis} and
 the very recent papers of Ben Alaya and Kebaier \cite{BenKeb1}, \cite{BenKeb2} (MLE).
In Ben Alaya and Kebaier \cite{BenKeb1}, \cite{BenKeb2} one can find a systematic study of the asymptotic behavior
 of the quadruplet \ $\big( \log(Y_t), Y_t, \int_0^t Y_s\,\dd s, \int_0^t 1/Y_s\,\dd s\big)$ \ as \ $t\to\infty$.
\ Finally, we note that Li and Ma \cite{LiMa} started to investigate the asymptotic behaviour of the (weighted) conditional
 LSE of the drift parameters for a CIR model driven by a stable noise (they call it a stable CIR model)
 from some discretely observed low frequency real data set.
To give another example besides the one-dimensional CIR process,
 we mention a model that is somewhat related to \eqref{2dim_affine} and parameter estimation
 of the appearing parameters based on continuous time observations has been considered.
It is the so-called Ornstein--Uhlenbeck process driven by \ $\alpha$-stable L\'evy motions, i.e.,
 \[
    \dd U_t = (m-\theta U_t)\,\dd t + \dd Z_t,\qquad t\geq 0,
 \]
 where \ $\theta>0$, $m\ne 0$, \ and \ $(Z_t)_{t\geq 0}$ \ is an \ $\alpha$-stable L\'evy motion with \ $\alpha\in(1,2)$.
\ For this model Hu and Long investigated the question of parameter estimation, see \cite{HuLon1}, \cite{HuLon2} and \cite{HuLon3}.

It would be possible to calculate the discretized version of the estimators
 presented in this paper using the same procedure as in Ben Alaya and Kebaier
 \cite[Section 4]{BenKeb1} valid for discrete time observations of high
 frequency.
However, it is out of the scope of the present paper.

We give a brief overview of the structure of the paper.
Section \ref{Section_preliminaires} is devoted to a preliminary discussion of
 the existence and uniqueness of a strong solution of the SDE \eqref{2dim_affine},
 we make a classification of the model (see Definition \ref{Def_criticality}),
 we also recall our results in Barczy et al.~\cite{BarDorLiPap2} on
 the existence of a unique stationary distribution and ergodicity
 for the affine process given by SDE \eqref{2dim_affine}, see Theorem \ref{Thm_ergodic12}.
Further, we recall some limit theorems for
 continuous local martingales that will be used later on for studying the asymptotic behaviour of
 the MLE of \ $(a,b,m,\theta)$ \ and the LSE of \ $(m,\theta)$, \ respectively.
\ In Sections \ref{section_EUMLE}--\ref{section_ALSE}
we study the asymptotic behavior of the MLE of \ $(a,b,m,\theta)$ \ and LSE of  \ $(m,\theta)$ \ proving that
 the estimators are strongly consistent and asymptotically normal under appropriate conditions on the parameters.
We call the attention that for the MLE of \ $(a,b,m,\theta)$ \ we require a continuous time observation
 \ $(Y_t,X_t)_{t\in[0,T]}$ \ of the process \ $(Y,X)$, \  but for the LSE of \ $(m,\theta)$ \ we need
 a continuous time observation \ $(X_t)_{t\in[0,T]}$ \ only for the process \ $X$.
\ We note that in the critical case we obtained a different limit behaviour for the LSE of \ $(m,\theta)$,
 \ see Barczy et al. \cite[Theorem 3.2]{BarDorLiPap}.

A common generalization of the model \eqref{2dim_affine} and the well-known
 Heston model \cite{Hes} is a general affine diffusion two factor model
 \begin{align}\label{2dim_affine_gen}
   \begin{cases}
    \dd Y_t = (a-bY_t)\,\dd t + \sigma_1\sqrt{Y_t}\,\dd L_t,\\
    \dd X_t = (m-\kappa Y_t-\theta X_t)\,\dd t
              + \sigma_1\sqrt{Y_t}\,(\varrho L_t + \sqrt{1-\varrho^2}\dd B_t),
   \end{cases}
   t\geq 0,
 \end{align}
 where \ $a,\sigma_1,\sigma_2>0$, \ $b, m,\kappa, \theta \in \RR$,
 \ $\varrho\in(-1,1)$, \ and \ $(L_t)_{t\geq 0}$ \ and \ $(B_t)_{t\geq 0}$ \ are
 independent standard Wiener processes.
One does not need to estimate the parameters \ $\sigma_1$, \ $\sigma_2$ \ and
 \ $\varrho$, \ since these parameters could ---in principle, at least--- be
 determined (rather than estimated) using an arbitrarily short continuous time
 observation of \ $(Y,X)$, \ see Remark 2.5 in Barczy and Pap \cite{BarPap}.
For studying the parameter estimation of \ $a$, $b$, $m$, $\kappa$ \ and
 \ $\theta$ \ in the subcritical case, one needs to investigate ergodicity
 properties of the model \eqref{2dim_affine_gen}.
For the submodel \eqref{2dim_affine}, this has been proved in
 Barczy et al.~\cite{BarDorLiPap2}, see also Theorem \ref{Thm_ergodic12}.
For the Heston model, ergodicity of the first coordinate process \ $Y$ \ is
 sufficient for statistical purposes, see Barczy and Pap \cite{BarPap}; the
 existence of a unique stationary distribution and the ergodicity for \ $Y$
 \ has been proved by Cox et al.~\cite[Equation (20)]{CoxIngRos} and Li and Ma
 \cite[Theorem 2.6]{LiMa}.
After showing appropriate ergodicity properties of the model
 \eqref{2dim_affine_gen}, one could obtain the asymptotic behavior of the MLE
 and LSE of \ $(a,b,m,\kappa,\theta)$ \ with a similar method used in the
 present paper.

\section{Preliminaires}\label{Section_preliminaires}

Let \ $\NN$, \ $\ZZ_+$, \ $\RR$ \ and \ $\RR_+$ \ denote the sets of positive integers, non-negative integers, real numbers and
 non-negative real numbers, respectively.
\ By \ $\|x\|$ \ and \ $\|A\|$ \ we denote the Euclidean norm of a vector \ $x\in\RR^m$ \ and the induced matrix norm
 \ $\Vert A\Vert=\sup\{\Vert Ax\Vert : x\in\RR^m, \; \Vert x\Vert=1\}$ \ of a matrix \ $A\in\RR^{n\times m}$, \ respectively.

The next proposition is about the existence and uniqueness of a strong solution of the SDE \eqref{2dim_affine}
 which follows from the theorem due to Yamada and Watanabe, further it clarifies that the process given by the SDE \eqref{2dim_affine}
 belongs to the family of regular affine processes, see Barczy et al.~\cite[Theorem 2.2 with \ $\alpha=2$]{BarDorLiPap2}.

\begin{Pro}\label{Pro_affine}
Let \ $\bigl(\Omega, \cF, \PP\bigr)$ \ be a probability space, and
 \ $(L_t, B_t)_{t\in\RR_+}$ \ be a 2-dimensional standard Wiener process.
Let \ $(\eta_0,\zeta_0)$ \ be a random vector independent of \ $(L_t,B_t)_{t\geq 0}$ \ satisfying
 \ $\PP(\eta_0\geq 0) = 1$.
\ Then, for all \ $a>0$, \ and \ $b, m, \theta\in \RR$, \ there is a (pathwise)
 unique strong solution \ $(Y_t,X_t)_{t\geq 0}$ \ of the SDE \eqref{2dim_affine} such that \ $\PP((Y_0,X_0) = (\eta_0,\zeta_0))=1$ \ and
 \ $\PP(Y_t\geq 0,\; \forall\;t\geq 0)=1$.
\ Further, we have
  \begin{align}\label{help50}
   Y_t = \ee^{-b(t-s)}\left( Y_s + a\int_s^t \ee^{-b(s-u)}\,\dd u
                                 + \int_s^t \ee^{-b(s-u)}\sqrt{Y_u} \,\dd L_u\right),
                                 \qquad 0\leq s\leq t,
 \end{align}
 and
 \begin{align}\label{help2}
  X_t = \ee^{-\theta(t-s)}\left( X_s + m\int_s^t \ee^{-\theta(s-u)}\,\dd u
                                 + \int_s^t \ee^{-\theta(s-u)}\sqrt{Y_u} \,\dd B_u\right),
                                 \qquad 0\leq s\leq t.
 \end{align}
Moreover, \ $(Y_t,X_t)_{t\geq 0}$ \ is a regular affine process with infinitesimal generator
 \begin{align*}
  (\cA f)(y,x)
   = (a-by)f_1'(y,x) + (m-\theta x)f_2'(y,x)
     + \frac{1}{2}y(f_{1,1}''(y,x) + f_{2,2}''(y,x)),
 \end{align*}
 where \ $(y,x)\in\RR_+\times\RR$, \ $f\in\cC^2_c(\RR_+\times\RR, \RR)$,
 \ $f_i'$, \ $i=1,2$, \ and \ $f_{i,j}''$, \ $i,j\in\{1,2\}$,
 \ denote the first and second order partial derivatives of \ $f$ \ with respect to
 its \ $i$-th and \ $i$-th and \ $j$-th variables, respectively, and \ $\cC^2_c(\RR_+\times\RR, \RR)$ \
 is the set of twice continuously differentiable real-valued functions defined on \ $\RR_+\times\RR$ \
 having compact support.
\end{Pro}

\begin{Rem}\label{augment}
Note that in Proposition \ref{Pro_affine} the unique strong solution \ $(Y_t,X_t)_{t\geq 0}$ \ of the SDE
 \eqref{2dim_affine} is adapted to the augmented filtration
 \ $(\cF_t)_{t\in\RR_+}$ \ corresponding to \ $(L_t, B_t)_{t\in\RR_+}$ \ and
 \ $(\eta_0,\zeta_0)$, \ constructed as in Karatzas and Shreve
 \cite[Section 5.2]{KarShr}.
Note also that \ $(\cF_t)_{t\in\RR_+}$ \ satisfies the usual conditions, i.e.,
 the filtration \ $(\cF_t)_{t\in\RR_+}$
 \ is right-continuous and \ $\cF_0$ \ contains all the $\PP$-null sets in \ $\cF$.
\ Further, \ $(L_t)_{t\in\RR_+}$ \ and \ $(B_t)_{t\in\RR_+}$ \ are independent \ $(\cF_t)_{t\in\RR_+}$-standard Wiener processes.
In Proposition \ref{Pro_affine} it is the assumption \ $a>0$ \ which
 ensures \ $\PP(Y_t\geq 0,\; \forall\;t\geq 0)=1$.
\proofend
\end{Rem}

In what follows we will make a classification of the affine processes given by the SDE \eqref{2dim_affine}.
\ First we recall a result about the first moment of \ $(Y_t,X_t)_{t\in\RR_+}$, \ see Proposition 3.2 in Barczy et al.
 \cite{BarDorLiPap}.

\begin{Pro}\label{Pro_moments}
Let \ $(Y_t,X_t)_{t\in\RR_+}$ \ be an affine diffusion process given by the SDE \eqref{2dim_affine} with a random initial value
 \ $(\eta_0,\zeta_0)$ \ independent of \ $(L_t,B_t)_{t\geq 0}$ \ such that \ $\PP(\eta_0\geq 0) = 1$,
 \ $\EE(\eta_0)<\infty$ \ and \ $\EE(\vert\zeta_0\vert)<\infty$.
\ Then
 \begin{align*}
  \begin{bmatrix}
   \EE(Y_t) \\
   \EE(X_t) \\
  \end{bmatrix}
  = \begin{bmatrix}
     \ee^{-bt} & 0 \\
     0 & \ee^{-\theta t} \\
    \end{bmatrix}
    \begin{bmatrix}
     \EE(\eta_0) \\
     \EE(\zeta_0) \\
    \end{bmatrix}
    + \begin{bmatrix}
       \int_0^t \ee^{-bs} \, \dd s & 0 \\
       0 & \int_0^t \ee^{-\theta s} \, \dd s \\
      \end{bmatrix}
      \begin{bmatrix}
       a \\
       m \\
      \end{bmatrix},\qquad t\in\RR_+.
 \end{align*}
\end{Pro}

Proposition \ref{Pro_moments} shows that the asymptotic behavior of the first moment of \ $(Y_t,X_t)_{t\in\RR_+}$ \ as \ $t\to\infty$ \
 is determined by the spectral radius of the diagonal matrix \ $\diag(\ee^{-bt}, \ee^{-\theta t})$, \ which motivates
 our classification of the affine processes given by the SDE \eqref{2dim_affine}.

\begin{Def}\label{Def_criticality}
Let \ $(Y_t,X_t)_{t\in\RR_+}$ \ be an affine diffusion process given by the SDE \eqref{2dim_affine} with a random initial value
 \ $(\eta_0,\zeta_0)$ \ independent of \ $(L_t,B_t)_{t\geq 0}$ \ satisfying
 \ $\PP(\eta_0\geq 0) = 1$.
\ We call \ $(Y_t,X_t)_{t\in\RR_+}$ \ subcritical, critical or supercritical if the spectral radius of the matrix
 \[
   \begin{bmatrix}
    \ee^{-bt} & 0 \\
    0 & \ee^{-\theta t} \\
   \end{bmatrix}
 \]
 is less than \ $1$, \ equal to \ $1$ \ or greater than \ $1$, \ respectively.
\end{Def}

Note that, since the spectral radius of the matrix given in Definition \ref{Def_criticality}
 is \ $\max(\ee^{-bt},\ee^{-\theta t})$, \ the affine process given in Definition \ref{Def_criticality}
 is
 \begin{align*}
  \text{subcritical} \qquad & \text {if \ $b>0$ \ and \ $\theta>0$,}\\
  \text{critical} \qquad
  & \text{if \ $b=0$, \ $\theta\geq 0$ \ or \ $b\geq 0$, \ $\theta=0$,}\\
  \text{supercritical} \qquad & \text{if \ $b<0$ \ or \ $\theta<0$.}
 \end{align*}
Further, under the conditions of Proposition \ref{Pro_moments}, by an easy calculation,
 if \ $b>0$ \ and \ $\theta>0$, \ then
 \[
  \lim_{t\to\infty}
  \begin{bmatrix}
   \EE(Y_t) \\
   \EE(X_t) \\
  \end{bmatrix}
   = \begin{bmatrix}
      \frac{a}{b} \\
      \frac{m}{\theta} \\
  \end{bmatrix},
 \]
 if \ $b=0$ \ and \ $\theta=0$, \ then
 \[
  \lim_{t\to\infty}
  \begin{bmatrix}
   \frac{1}{t}\EE(Y_t) \\
   \frac{1}{t}\EE(X_t) \\
  \end{bmatrix}
   = \begin{bmatrix}
      a \\
      m \\
  \end{bmatrix},
 \]
 if \ $b=0$ \ and \ $\theta>0$, \ then
 \[
  \lim_{t\to\infty}
  \begin{bmatrix}
   \frac{1}{t}\EE(Y_t) \\
              \EE(X_t) \\
  \end{bmatrix}
   = \begin{bmatrix}
      a \\
      \frac{m}{\theta} \\
  \end{bmatrix},
 \]
 if \ $b>0$ \ and \ $\theta=0$, \ then
 \[
  \lim_{t\to\infty}
  \begin{bmatrix}
    \EE(Y_t) \\
   \frac{1}{t}\EE(X_t) \\
  \end{bmatrix}
   = \begin{bmatrix}
      \frac{a}{b} \\
      m \\
  \end{bmatrix},
 \]
 and if \ $b<0$ \ and \ $\theta<0$, \ then
 \[
  \lim_{t\to\infty}
  \begin{bmatrix}
   \ee^{bt}\EE(Y_t) \\
   \ee^{\theta t}\EE(X_t) \\
  \end{bmatrix}
   = \begin{bmatrix}
      \EE(\eta_0) - \frac{a}{b} \\
      \EE(\zeta_0) - \frac{m}{\theta} \\
   \end{bmatrix}.
 \]
Remark also that Definition \ref{Def_criticality} of criticality is in accordance with the
 corresponding definition for one-dimensional continuous state branching processes, see, e.g.,
 Li \cite[page 58]{Li}.

In the sequel \ $\stoch$ \ and \ $\distr$ \ will denote convergence in probability and in distribution, respectively.

The following result states the existence of a unique stationary distribution and the ergodicity
 for the affine process given by the SDE \eqref{2dim_affine}, see Theorems 3.1 with \ $\alpha=2$ \
 and Theorem 4.2 in Barczy et al.~\cite{BarDorLiPap2}.

\begin{Thm}\label{Thm_ergodic12}
Let us consider the 2-dimensional affine model \eqref{2dim_affine} with
 \ $a>0$, \ $b>0$, $m\in\RR$, \ $\theta>0$, \ and with a random initial value
 \ $(\eta_0,\zeta_0)$ \ independent of \ $(L_t,B_t)_{t\geq 0}$ \ satisfying \ $\PP(\eta_0\geq 0) = 1$.
\ Then
 \renewcommand{\labelenumi}{{\rm(\roman{enumi})}}
 \begin{enumerate}
 \item \ $(Y_t,X_t)\distr (Y_\infty,X_\infty)$ \ as \ $t\to\infty$, \ and the distribution of
       \ $(Y_\infty,X_\infty)$ \ is given by
       \begin{align}\label{help18}
          \EE\big(\ee^{-\lambda_1 Y_\infty + i\lambda_2 X_\infty}\big)
            = \exp\left\{-a\int_0^\infty v_s(\lambda_1,\lambda_2)\,\dd s + i\frac{m}{\theta}\lambda_2\right\},
          \qquad (\lambda_1,\lambda_2)\in\RR_+\times\RR,
       \end{align}
       where \ $v_t(\lambda_1,\lambda_2)$, $t\geq 0$, \ is the unique
       non-negative solution of the (deterministic) differential
       equation
       \begin{align}\label{DE1}
         \begin{cases}
          \frac{\partial v_t}{\partial t} (\lambda_1,\lambda_2)
               = -bv_t(\lambda_1,\lambda_2) - \frac{1}{2}(v_t(\lambda_1,\lambda_2))^2
                  + \frac{1}{2}\ee^{-2\theta t}\lambda_2^2,\qquad t\geq 0,\\
          v_0(\lambda_1,\lambda_2) = \lambda_1.
         \end{cases}
       \end{align}
 \item supposing that the random initial value \ $(\eta_0,\zeta_0)$ \ has the same distribution as
       \ $(Y_\infty,X_\infty)$ \ given in part (i), we have \ $(Y_t,X_t)_{t\geq 0}$ \ is strictly stationary.
 \item for all Borel measurable functions \ $f:\RR^2\to\RR$ \ such that
       \ $\EE\big(\vert f(Y_\infty,X_\infty)\vert\big)<\infty$, \ we have
       \begin{equation}\label{help_ergodic}
         \PP\left( \lim_{T\to\infty} \frac{1}{T}\int_0^T f(Y_s,X_s)\,\dd s
           = \EE (f(Y_\infty,X_\infty)) \right)=1,
       \end{equation}
       where the distribution of \ $(Y_\infty,X_\infty)$ \ is given by \eqref{help18} and \eqref{DE1}.
 \end{enumerate}
 Moreover, the random variable \ $(Y_\infty,X_\infty)$ \ is absolutely continuous,
 the Laplace transform of \ $Y_\infty$ \ takes the form
 \begin{align}\label{help74}
  \EE(\ee^{-\lambda_1 Y_\infty}) = \left(1 + \frac{\lambda_1}{2b}\right)^{-2a},
           \qquad \lambda_1\in\RR_+,
 \end{align}
 i.e., \ $Y_\infty$ \ has Gamma distribution with parameters \ $2a$ \ and
 \ $2b$, \ all the (mixed) moments of \ $(Y_\infty,X_\infty)$ \ of any order
 are finite, i.e., \ $\EE(Y_\infty^n\vert X_\infty\vert^p)<\infty$
 \ for all \ $n,p\in\ZZ_+$, \ and especially,
 \begin{align*}
          &\EE(Y_\infty)=\frac{a}{b}, \qquad \EE(X_\infty)=\frac{m}{\theta}, \\
          &\EE(Y_\infty^2) = \frac{a(2a+1)}{2b^2},
         \qquad
         \EE(Y_\infty X_\infty) = \frac{ma}{\theta b},
         \qquad
         \EE(X_\infty^2) = \frac{a\theta+ 2bm^2}{2b\theta^2},\\
         &\EE(Y_\infty X_\infty^2)
            = \frac{a}{(b+2\theta)2b^2\theta^2}
               \big(\theta(ab+2a\theta+\theta) + 2m^2b(2\theta + b)\big).
 \end{align*}
\end{Thm}

In all what follows we will suppose that we have continuous time observations for the process \ $(Y,X)$,
 \ i.e., \ $(Y_t,X_t)_{t\in[0,T]}$ \ can be observed for some \ $T>0$,
 \ and our aim is to deal with parameter estimation of \ $(a,b,m,\theta)$.
\ We also deal with parameter estimation of \ $\theta$ \ provided that the parameter \ $m\in\RR$ \ is supposed to be known.

Next we recall some limit theorems for continuous local martingales.
We will use these limit theorems in the sequel for studying the asymptotic behaviour of different kinds of estimators
 for \ $(a,b,m,\theta)$.
\ First we recall a strong law of large numbers for continuous local martingales, see, e.g.,
 Liptser and Shiryaev \cite[Lemma 17.4]{LipShiII}.

\begin{Thm}\label{DDS_stoch_int}
Let \ $\big(\Omega, \cF, (\cF_t)_{t\geq 0}, \PP\big)$ \ be a filtered probability space satisfying the usual conditions.
Let \ $(M_t)_{t \geq 0}$ \ be a square-integrable continuous local martingale
 with respect to the filtration \ $(\cF_t)_{t\geq 0}$ \ started from \ $0$.
\ Let \ $(\xi_t)_{t\geq 0}$ \ be a progressively measurable process such that
 \[
   \PP\left(\int_0^t(\xi_u)^2\,\dd\langle M\rangle_u<\infty\right)=1,
   \quad t\geq 0,
 \]
 and
 \begin{align}\label{SEGED_STRONG_CONSISTENCY2}
  \PP\left(\lim_{t\to\infty}\int_0^t(\xi_u)^2\,\dd\langle M\rangle_u=\infty\right)
  =1,
 \end{align}
 where \ $(\langle M\rangle_t)_{t\geq 0}$ \ denotes the quadratic variation process of \ $M$.
\ Then
 \begin{align}\label{SEGED_STOCH_INT_SLLN}
   \PP\left(\lim_{t\to\infty}\frac{\int_0^t\xi_u\,\dd M_u}
         {\int_0^t(\xi_u)^2\,\dd \langle M\rangle_u}=0\right)=1.
 \end{align}
In case of \ $M_t=B_t$, \ $t\geq 0$, \ where \ $(B_t)_{t\geq 0}$ \ is a
 standard Wiener process, the progressive measurability of \ $(\xi_t)_{t\geq 0}$ \ can be relaxed
 to measurability and adaptedness to the filtration \ $(\cF_t)_{t\geq 0}$.
\end{Thm}

The next theorem is about the asymptotic behaviour of continuous multivariate local martingales.

\begin{Thm}{\bf (van Zanten \cite[Theorem 4.1]{Zan})}\label{THM_Zanten}
 Let \ $\big(\Omega, \cF, (\cF_t)_{t\geq 0}, \PP\big)$ \ be a filtered probability
 space satisfying the usual conditions.
Let \ $(M_t)_{t \geq 0}$ \ be a \ $d$-dimensional square-integrable continuous local martingale with respect to the
 filtration \ $(\cF_t)_{t\geq 0}$ \ started from \ $0$.
\ Suppose that there exists a function \ $Q:[0,\infty)\to \RR^{d\times d}$ \ such
 that \ $Q(t)$ \ is a non-random, invertible matrix for all \ $t\geq 0$,
 \ $\lim_{t\to\infty}\Vert Q(t)\Vert=0$ \ and
 \[
  Q(t)\langle M\rangle_tQ(t)^\top\stoch \eta\eta^\top
    \quad \text{as} \quad  t\to\infty,
 \]
 where \ $\eta$ \ is a \ $d\times d$ \ random matrix defined on
 \ $\big(\Omega, \cF, \PP\big)$.
\ Then, for each \ $\RR^k$-valued random variable \ $V$ \ defined on
 \ $\big(\Omega, \cF, \PP\big)$, \ it holds that
 \[
       (Q(t)M_t,V)\distr (\eta Z,V) \qquad\text{as}\qquad t\to\infty,
 \]
 where \ $Z$ \ is a \ $d$-dimensional standard normally distributed random variable independent
 of \ $(\eta,V)$.
\end{Thm}

We note that Theorem \ref{THM_Zanten} remains true if the function \ $Q$, \ instead of the interval \ $[0,\infty)$,
 \ is defined only on an interval \ $[t_0,\infty)$ \ with some \ $t_0>0$.

\section{Existence and uniqueness of maximum likelihood estimator}\label{section_EUMLE}

Let \ $\PP_{(a,b,m,\theta)}$ \ denote the probability measure on the measurable space
 \ $(\cC(\RR_+, \RR_+\times\RR),\cB(\cC(\RR_+, \RR_+\times\RR)))$ \ induced by the process \ $(Y_t,X_t)_{t\geq 0}$ \
 corresponding to the parameters \ $(a,b,m,\theta)$ \ and initial value \ $(Y_0,X_0)$.
\ Here \ $\cC(\RR_+, \RR_+\times\RR)$ \ denotes the set of continuous \ $\RR_+\times\RR$-valued functions defined
 on \ $\RR_+$, \ $\cB(\cC(\RR_+, \RR_+\times\RR))$ \ is the Borel \ $\sigma$-algebra on it, and
 we suppose that the space \ $(\cC(\RR_+, \RR_+\times\RR),\cB(\cC(\RR_+, \RR_+\times\RR)))$ \ is endowed with
 the natural filtration \ $(\cA_t)_{t\geq 0}$, \ given by \ $\cA_t:=\varphi_t^{-1}(\cB(\cC(\RR_+, \RR_+\times\RR)))$, \ where
 \ $\varphi_t:\cC(\RR_+, \RR_+\times\RR)\to \cC(\RR_+, \RR_+\times\RR)$ \ is the mapping \ $\varphi_t(f)(s):=f(t\wedge s)$, \ $s\geq 0$.
\ For all \ $T>0$, \ let \ $\PP_{(a,b,m,\theta),T}:=\PP_{(a,b,m,\theta)}\vert_{\cA_T}$ \ be the restriction of
 \ $\PP_{(a,b,m,\theta)}$ \ to \ $\cA_T$.

\begin{Lem}\label{LEMMA_MLE_likelihood}
Let \ $a\geq 1/2$, \ $b,m,\theta\in\RR$, \ $T>0$, \ and suppose that \ $\PP(Y_0>0)=1$.
\ Let \ $\PP_{(a,b,m,\theta)}$ \ and \ $\PP_{(1,0,0,0)}$ \ denote the probability measures
 induced by the unique strong solutions of the SDE \eqref{2dim_affine}
 corresponding to the parameters \ $(a,b,m,\theta)$ \ and \ $(1,0,0,0)$ \ with the same initial value
 \ $(Y_0,X_0)$, \ respectively.
\ Then \ $\PP_{(a,b,m,\theta),T}$ \ and \ $\PP_{(1,0,0,0),T}$ \ are absolutely continuous with respect to
 each other, and the Radon-Nykodim derivative of \ $\PP_{(a,b,m,\theta),T}$ \ with respect to
 \ $\PP_{(1,0,0,0),T}$ \ (so called likelihood ratio) takes the form
 \begin{align*}
  L^{(a,b,m,\theta),(1,0,0,0)}_T((Y_s,X_s)_{s\in[0,T]})
  =\exp\Bigg\{&\int_0^T
                \left(\frac{a-bY_s-1}{Y_s} \, \dd Y_s
                      + \frac{m-\theta X_s}{Y_s} \, \dd X_s\right) \\
              &-\frac{1}{2}
                \int_0^T
                 \frac{(a-bY_s-1)(a-bY_s+1) + (m-\theta X_s)^2}{Y_s} \, \dd s
       \Bigg\} ,
 \end{align*}
 where \ $(Y_t,X_t)_{t\geq 0}$ \ denotes the unique strong solution of the SDE \eqref{2dim_affine}
 corresponding to the parameters \ $(a,b,m,\theta)$ \ and the initial value \ $(Y_0,X_0)$.
\end{Lem}

\noindent{\bf Proof.}
First note that the SDE \eqref{2dim_affine} can be written in the form:
 \begin{align*}
  \dd \begin{bmatrix}
        Y_t \\
        X_t \\
      \end{bmatrix}
   = \left[
       \begin{bmatrix}
         -b & 0 \\
         0 & -\theta \\
       \end{bmatrix}
       \begin{bmatrix}
        Y_t \\
        X_t \\
      \end{bmatrix}
      + \begin{bmatrix}
        a \\
        m \\
      \end{bmatrix}
      \right]\,\dd t
      + \begin{bmatrix}
         \sqrt{Y_t} & 0 \\
         0 & \sqrt{Y_t} \\
       \end{bmatrix}
       \begin{bmatrix}
        \dd L_t \\
        \dd B_t \\
      \end{bmatrix},
     \qquad t\geq 0.
 \end{align*}
Note also that under the condition \ $a\geq\frac{1}{2}$, \ for all \ $y_0>0$,
 \ we have \ $\PP(Y_t>0,\, \forall\;t\in\RR_+ \mid Y_0 =y_0)=1$, \ see, e.g.,
 page 442 in Revuz and Yor \cite{RevYor}.
Since \ $\PP(Y_0>0)=1$, \ by the law of total probability,
 \ $\PP(Y_t>0,\,\forall\, t\in\RR_+)=1$.

We intend to use Lemma \ref{LEMMA_LipShi}.
By Proposition \ref{Pro_affine}, under the conditions of the present lemma,
 there is a pathwise unique strong solution of the SDE \eqref{2dim_affine}.
We has to check
 \begin{align*}
   \int_0^T \frac{(a-bY_s)^2 + (m-\theta X_s)^2 + 1}{Y_s} \,\dd s<\infty
                  \qquad \text{a.s.~for all \ $T\in\RR_+$.}
 \end{align*}
Since \ $(Y,X)$ \ has continuous sample paths almost surely, this holds if
 \begin{align}\label{help78}
    \int_0^T \frac{1}{Y_s}\,\dd s <\infty \qquad \text{a.s.~for all \ $T\in\RR_+$.}
 \end{align}
Under the conditions \ $a\geq 1/2$ \ and \ $\PP(Y_0>0)=1$,
 \ Theorems 1 and 3 in Ben-Alaya and Kebaier \cite{BenKeb2} yield \eqref{help78}.
More precisely, if \ $a\geq\frac{1}{2}$ \ and \ $y_0>0$, \ then Theorems 1 and 3 in Ben-Alaya and Kebaier \cite{BenKeb2} yield
 \[
    \PP\left(\int_0^T \frac{1}{Y_s} \,\dd s <\infty \,\Big\vert\, Y_0 = y_0\right)=1,\qquad T\in\RR_+.
 \]
Since \ $\PP(Y_0>0)=1$, \ by the law of total probability, we get \eqref{help78}.
We give another, direct proof for \eqref{help78}.
Namely, since \ $Y$ \ has continuous sample paths almost surely and \ $\PP(Y_t>0, \, \forall\,t\in\RR_+)=1$, \
 we have \ $\PP(\inf_{t\in[0,T]} Y_t >0)=1$ \ for all \ $T\in\RR_+$, \ which yields \eqref{help78}.
\proofend

By Lemma \ref{LEMMA_MLE_likelihood}, under its conditions the log-likelihood function takes the form
  \begin{align*}
   \log L^{(a,b,m,\theta),(1,0,0,0)}_T((Y_s,X_s)_{s\in[0,T]})
     &= (a-1) \int_0^T \frac{1}{Y_s}\,\dd Y_s - b (Y_T-Y_0)
        + m \int_0^T \frac{1}{Y_s}\,\dd X_s\\
     &\phantom{=\;}
        - \theta \int_0^T \frac{X_s}{Y_s}\,\dd X_s
        - \frac{a^2-1}{2}\int_0^T \frac{1}{Y_s}\,\dd s
        + abT
       - \frac{b^2}{2}\int_0^T Y_s\,\dd s\\
     &\phantom{=\;}
       - \frac{m^2}{2}\int_0^T \frac{1}{Y_s}\,\dd s
       + m\theta\int_0^T \frac{X_s}{Y_s}\,\dd s
       -\frac{\theta^2}{2} \int_0^T \frac{X_s^2}{Y_s}\,\dd s\\
     &=:f_T(a,b,m,\theta),
   \qquad T>0.
 \end{align*}

We remark that for all \ $T>0$ \ and all initial values \ $(Y_0,X_0)$, \ the
 probability measures \ $\PP_{(a,b,m,\theta),T}$,\ $a\geq 1/2$, $b,m,\theta\in\RR$, \ are absolutely continuous with respect to each other,
 and hence it does not matter which measure is taken as a reference measure for defining the MLE
 (we have chosen \ $\PP_{(1,0,0,0),T}$).
\ For more details, see, e.g., Liptser and Shiryaev \cite[page 35]{LipShiI}.
Then the equation \ $\frac{\partial f_T}{\partial \theta}(a,b,m,\theta) = 0$
 \ takes the form
  \begin{align*}
   -\int_0^T \frac{X_s}{Y_s}\,\dd X_s + m\int_0^T \frac{X_s}{Y_s}\,\dd s
   - \theta \int_0^T \frac{X_s^2}{Y_s}\,\dd s
   = 0 .
  \end{align*}
Moreover, the system of equations
 \begin{align*}
  \frac{\partial f_T}{\partial a}(a,b,m,\theta) = 0,
  \qquad
  \frac{\partial f_T}{\partial b}(a,b,m,\theta) = 0,
  \qquad
  \frac{\partial f_T}{\partial m}(a,b,m,\theta) = 0,
  \qquad
  \frac{\partial f_T}{\partial \theta}(a,b,m,\theta) = 0,
 \end{align*}
 takes the form
 \begin{align*}
  \begin{bmatrix}
   \int_0^T \frac{1}{Y_s}\,\dd s & -T & 0 & 0 \\
   -T &  \int_0^T Y_s\,\dd s & 0 & 0 \\
   0 & 0 & \int_0^T \frac{1}{Y_s}\,\dd s & -\int_0^T \frac{X_s}{Y_s}\,\dd s \\
   0 & 0 & -\int_0^T \frac{X_s}{Y_s}\,\dd s & \int_0^T \frac{X_s^2}{Y_s}\,\dd s
  \end{bmatrix}
  \begin{bmatrix}
   a \\
   b \\
   m \\
   \theta
  \end{bmatrix}
  =
  \begin{bmatrix}
   \int_0^T \frac{1}{Y_s}\,\dd Y_s \\
   -(Y_T-Y_0) \\
   \int_0^T \frac{1}{Y_s}\,\dd X_s \\
   -\int_0^T \frac{X_s}{Y_s}\,\dd X_s
  \end{bmatrix}.
 \end{align*}

First, we suppose that \ $a\geq 1/2$, \ and \ $b\in\RR$ \ and \ $m\in\RR$
 \ are known.
By maximizing \ $\log L^{(a,b,m,\theta),(1,0,0,0)}_T$ \ in \ $\theta\in\RR$, \ we get the MLE
 of \ $\theta$ \ based on the observations \ $(Y_t,X_t)_{t\in[0,T]}$,
 \begin{align}\label{MLE_theta}
  {\widetilde\theta}^{\mathrm{MLE}}_T
  := \frac{-\int_0^T \frac{X_s}{Y_s}\,\dd X_s + m\int_0^T \frac{X_s}{Y_s}\,\dd s}
          {\int_0^T \frac{X_s^2}{Y_s}\,\dd s},\qquad T>0,
 \end{align}
 provided that \ $\int_0^T \frac{X_s^2}{Y_s}\,\dd s>0$.
\ Indeed,
 \[
   \frac{\partial^2 f_T}{\partial \theta^2}(\theta,m)
   = - \int_0^T \frac{X_s^2}{Y_s}\,\dd s
   < 0 .
 \]
Using the SDE \eqref{2dim_affine}, one can also get
 \begin{align}\label{help52}
  {\widetilde\theta}^{\mathrm{MLE}}_T  - \theta
  = - \frac{\int_0^T \frac{X_s}{\sqrt{Y_s}}\,\dd B_s}
           {\int_0^T \frac{X_s^2}{Y_s}\,\dd s},
  \qquad T > 0 ,
 \end{align}
 provided that \ $\int_0^T \frac{X_s^2}{Y_s}\,\dd s>0$.
\ Note that the estimator \ ${\widetilde\theta}^{\mathrm{MLE}}_T$ \ does not
 depend on the parameters \ $a\geq 1/2$ \ and \ $b\in\RR$.
\ In fact, if we maximize \ $\log L^{(a,b,m,\theta),(1,0,0,0)}_T$ \ in
 \ $(a,b,\theta)\in\RR^3$, \ then we obtain the MLE of \ $(a,b,\theta)$
 \ supposing that \ $m \in \RR$ \ is known, and one can observe that the MLE
 of \ $\theta$ \ by this procedure coincides with \ ${\widetilde\theta}^{\mathrm{MLE}}_T$.

By maximizing \ $\log L^{(a,b,m,\theta),(1,0,0,0)}_T$ \ in \ $(a,b,m,\theta)\in\RR^{4}$,
 \ the MLE of \ $(a,b,m,\theta)$ \ based on the observations \ $(Y_t,X_t)_{t\in[0,T]}$ \ takes the form
 \begin{align}\label{MLE_theta_m_1}
  \widehat{a}_T^{\mathrm{MLE}}
  &:=\frac{\int_0^T Y_s\,\dd s \int_0^T \frac{1}{Y_s}\,\dd Y_s - T(Y_T-Y_0)}
          {\int_0^T Y_s\,\dd s \int_0^T \frac{1}{Y_s}\,\dd s - T^2},\qquad T>0,\\\label{MLE_theta_m_2}
  \widehat{b}_T^{\mathrm{MLE}}
  &:=\frac{T\int_0^T \frac{1}{Y_s}\,\dd Y_s
           -(Y_T-Y_0) \int_0^T \frac{1}{Y_s}\,\dd s}
          {\int_0^T Y_s\,\dd s \int_0^T \frac{1}{Y_s}\,\dd s - T^2},\qquad T>0,\\\label{MLE_theta_m_3}
  \widehat m_T^{\mathrm{MLE}}
  &:=\frac{\int_0^T\frac{X_s^2}{Y_s}\,\dd s\int_0^T \frac{1}{Y_s}\,\dd X_s
           -\int_0^T \frac{X_s}{Y_s}\,\dd s \int_0^T \frac{X_s}{Y_s}\,\dd X_s}
          {\int_0^T \frac{X_s^2}{Y_s}\,\dd s \int_0^T \frac{1}{Y_s}\,\dd s - \left(\int_0^T \frac{X_s}{Y_s}\,\dd s\right)^2},\qquad T>0,\\
  \widehat\theta_T^{\mathrm{MLE}}
  &:=\frac{\int_0^T\frac{X_s}{Y_s}\,\dd s\int_0^T \frac{1}{Y_s}\,\dd X_s
           -\int_0^T \frac{1}{Y_s}\,\dd s \int_0^T \frac{X_s}{Y_s}\,\dd X_s}
          {\int_0^T \frac{X_s^2}{Y_s}\,\dd s \int_0^T \frac{1}{Y_s}\,\dd s
           - \left(\int_0^T \frac{X_s}{Y_s}\,\dd s\right)^2},\label{MLE_theta_m_4}
            \qquad T>0,
 \end{align}
 provided that \ $\int_0^T Y_s\,\dd s \int_0^T \frac{1}{Y_s}\,\dd s - T^2>0$ \ and
\ $\int_0^T \frac{X_s^2}{Y_s}\,\dd s \int_0^T \frac{1}{Y_s}\,\dd s - \left(\int_0^T \frac{X_s}{Y_s}\,\dd s\right)^2>0$.
\ Indeed,
 \begin{gather*}
  \begin{bmatrix}
   \frac{\partial^2 f_T}{\partial a^2}(a,b,m,\theta)
    & \frac{\partial^2 f_T}{\partial b \partial a}(a,b,m,\theta) \\
   \frac{\partial^2 f_T}{\partial a \partial b}(a,b,m,\theta)
    & \frac{\partial^2 f_T}{\partial b^2}(a,b,m,\theta)
  \end{bmatrix}
  =\begin{bmatrix}
    -\int_0^T \frac{1}{Y_s}\,\dd s & T \\
    T & -\int_0^T Y_s\,\dd s
   \end{bmatrix},\\
  \begin{bmatrix}
   \frac{\partial^2 f_T}{\partial m^2}(a,b,m,\theta)
    & \frac{\partial^2 f_T}{\partial \theta \partial m}(a,b,m,\theta) \\
   \frac{\partial^2 f_T}{\partial m \partial \theta}(a,b,m,\theta)
    & \frac{\partial^2 f_T}{\partial \theta^2}(a,b,m,\theta)
  \end{bmatrix}
  =\begin{bmatrix}
    -\int_0^T \frac{1}{Y_s}\,\dd s & \int_0^T \frac{X_s}{Y_s}\,\dd s \\
    \int_0^T \frac{X_s}{Y_s}\,\dd s & -\int_0^T \frac{X_s^2}{Y_s}\,\dd s
   \end{bmatrix},
 \end{gather*}
 and the positivity of \ $\int_0^T Y_s\,\dd s \int_0^T \frac{1}{Y_s}\,\dd s - T^2$ \ and
 \ $\int_0^T \frac{X_s^2}{Y_s}\,\dd s \int_0^T \frac{1}{Y_s}\,\dd s - \left(\int_0^T \frac{X_s}{Y_s}\,\dd s\right)^2$
 \ yield \ $\int_0^T \frac{1}{Y_s}\,\dd s >0$, \ respectively.
Using the SDE \eqref{2dim_affine} one can check that
 \begin{align*}
  \begin{bmatrix}
   \widehat{a}_T^{\mathrm{MLE}} - a \\
   \widehat{b}_T^{\mathrm{MLE}} - b \\
  \end{bmatrix}
  &=\begin{bmatrix}
     \int_0^T \frac{1}{Y_s}\,\dd s & -T \\
     -T & \int_0^T Y_s\,\dd s
    \end{bmatrix}^{-1}
    \begin{bmatrix}
     \int_0^T \frac{1}{\sqrt{Y_s}}\,\dd L_s \\
     -\int_0^T \sqrt{Y_s}\,\dd L_s
    \end{bmatrix},\\
  \begin{bmatrix}
   \widehat m_T^{\mathrm{MLE}} - m \\
   \widehat\theta_T^{\mathrm{MLE}} - \theta
  \end{bmatrix}
  &=\begin{bmatrix}
     \int_0^T \frac{1}{Y_s}\,\dd s & -\int_0^T \frac{X_s}{Y_s}\,\dd s \\
     -\int_0^T \frac{X_s}{Y_s}\,\dd s & \int_0^T \frac{X_s^2}{Y_s}\,\dd s
    \end{bmatrix}^{-1}
    \begin{bmatrix}
     \int_0^T \frac{1}{\sqrt{Y_s}}\,\dd B_s \\
     -\int_0^T \frac{X_s}{\sqrt{Y_s}}\,\dd B_s
    \end{bmatrix},
 \end{align*}
 and hence
 \begin{align}\label{help47a}
  \widehat{a}_T^{\mathrm{MLE}} - a
  &=\frac{\int_0^T Y_s\,\dd s \int_0^T \frac{1}{\sqrt{Y_s}}\,\dd L_s
          -T \int_0^T \sqrt{Y_s}\,\dd L_s}
         {\int_0^T Y_s\,\dd s \int_0^T \frac{1}{Y_s}\,\dd s - T^2},\qquad T>0,\\
  \label{help47b}
  \widehat{b}_T^{\mathrm{MLE}} - b
  &=\frac{T\int_0^T \frac{1}{\sqrt{Y_s}}\,\dd L_s
          -\int_0^T \frac{1}{Y_s}\,\dd s \int_0^T \sqrt{Y_s}\,\dd L_s}
         {\int_0^T Y_s\,\dd s \int_0^T \frac{1}{Y_s}\,\dd s - T^2},\qquad T>0,\\
  \label{help47}
  \widehat m_T^{\mathrm{MLE}} - m
  &=\frac{\int_0^T \frac{X_s^2}{Y_s}\,\dd s
          \int_0^T \frac{1}{\sqrt{Y_s}}\,\dd B_s
          -\int_0^T \frac{X_s}{Y_s}\,\dd s
           \int_0^T \frac{X_s}{\sqrt{Y_s}}\,\dd B_s}
         {\int_0^T \frac{X_s^2}{Y_s}\,\dd s \int_0^T \frac{1}{Y_s}\,\dd s
         -\left(\int_0^T \frac{X_s}{Y_s}\,\dd s\right)^2},
            \qquad T>0,
 \end{align}
 and
 \begin{align}\label{help48}
  \widehat\theta_T^{\mathrm{MLE}} - \theta
  =\frac{\int_0^T\frac{X_s}{Y_s}\,\dd s\int_0^T \frac{1}{\sqrt{Y_s}}\,\dd B_s
         -\int_0^T \frac{1}{Y_s}\,\dd s \int_0^T \frac{X_s}{\sqrt{Y_s}}\,\dd B_s}
        {\int_0^T \frac{X_s^2}{Y_s}\,\dd s \int_0^T \frac{1}{Y_s}\,\dd s
         -\left(\int_0^T \frac{X_s}{Y_s}\,\dd s\right)^2},\qquad T>0,
 \end{align}
 provided that \ $\int_0^T Y_s\,\dd s \int_0^T \frac{1}{Y_s}\,\dd s - T^2>0$ \ and \ $\int_0^T \frac{X_s^2}{Y_s}\,\dd s \int_0^T \frac{1}{Y_s}\,\dd s - \left(\int_0^T \frac{X_s}{Y_s}\,\dd s\right)^2>0$.

\begin{Rem}\label{Riemann}
For the stochastic integrals \ $\int_0^T \frac{1}{Y_s}\,\dd Y_s$, \ $\int_0^T \frac{X_s}{Y_s}\,\dd X_s$ \ and
 \ $\int_0^T \frac{1}{Y_s}\,\dd X_s$ \ in \eqref{MLE_theta_m_1}, \eqref{MLE_theta_m_2}, \eqref{MLE_theta_m_3} and
 \eqref{MLE_theta_m_4}, we have
 \begin{align}\label{measurability}
  \begin{split}
  & \sum_{i=1}^{\lfloor nT\rfloor } \frac{1}{Y_{\frac{i-1}{n}}}
             \big(Y_{\frac{i}{n}} - Y_{\frac{i-1}{n}} \big)
     \stoch \int_0^T \frac{1}{Y_s}\,\dd Y_s \qquad \text{as \ $n\to\infty$,}\\
  & \sum_{i=1}^{\lfloor nT\rfloor } \frac{X_{\frac{i-1}{n}}}{Y_{\frac{i-1}{n}}}
             \big(X_{\frac{i}{n}} - X_{\frac{i-1}{n}} \big)
     \stoch \int_0^T \frac{X_s}{Y_s}\,\dd X_s \qquad \text{as \ $n\to\infty$,}\\
  & \sum_{i=1}^{\lfloor nT\rfloor } \frac{1}{Y_{\frac{i-1}{n}}}
             \big(X_{\frac{i}{n}} - X_{\frac{i-1}{n}} \big)
    \stoch \int_0^T \frac{1}{Y_s}\,\dd X_s \qquad \text{as \ $n\to\infty$,}
   \end{split}
 \end{align}
 following from Proposition I.4.44 in Jacod and Shiryaev \cite{JSh} with the
 Riemann sequence of deterministic subdivisions
 \ $\left(\frac{i}{n} \land T\right)_{i\in\NN}$, \ $n \in \NN$.
\ Thus, there exist measurable functions \ $\Phi, \Psi, \Xi : C([0,T],\RR_+\times\RR) \to \RR$
 \ such that \ $\int_0^T \frac{1}{Y_s}\,\dd Y_s = \Phi((Y_s, X_s)_{s\in[0,T]})$,
 \ $\int_0^T \frac{X_s}{Y_s}\,\dd X_s = \Psi((Y_s, X_s)_{s\in[0,T]})$ \ and
 \ $\int_0^T \frac{1}{Y_s}\,\dd X_s = \Xi((Y_s, X_s)_{s\in[0,T]})$, \  since
 the convergences in \eqref{measurability} hold almost surely along suitable
 subsequences, the members of the sequences in
 \eqref{measurability} are measurable functions of \ $(Y_s, X_s)_{s\in[0,T]}$,
 \ and one can use Theorems 4.2.2 and 4.2.8 in Dudley \cite{Dud}.
Hence the right hand sides of \eqref{MLE_theta_m_1}, \eqref{MLE_theta_m_2},
 \eqref{MLE_theta_m_3} and \eqref{MLE_theta_m_4} are measurable functions of
 \ $(Y_s, X_s)_{s\in[0,T]}$, \ i.e., they are statistics.
\proofend
\end{Rem}

The next lemma is about the existence of \ $\widetilde \theta_T^{\mathrm{MLE}}$ \
 (supposing that \ $a\geq\frac{1}{2}$, \ $b\in\RR$ \ and \ $m\in\RR$ \ are known).

\begin{Lem}\label{LEMMA_MLE_exist_theta}
If \ $a\geq \frac{1}{2}$, \ $b,m,\theta\in\RR$, \ and \ $\PP(Y_0>0)=1$, \ then
 \begin{align}\label{help40}
  \PP\left( \int_0^T \frac{X_s^2}{Y_s}\,\dd s \in (0,\infty)\right)=1\quad
  \text{for all \ $T>0$,}
 \end{align}
 and hence there exists a unique MLE \ $\widetilde\theta^{\mathrm{MLE}}_T$ \ which has the form given in \eqref{MLE_theta}.
\end{Lem}

\noindent{\bf Proof.}
First note that under the condition \ $a\geq\frac{1}{2}$, \ for all \ $y_0>0$, \ we have
 \ $\PP(Y_t>0,\, \forall\;t\in\RR_+ \mid Y_0 =y_0)=1$, \ see, e.g., page 442 in Revuz and Yor \cite{RevYor}.
Since \ $\PP(Y_0>0)=1$, \ by the law of total probability, we get
 \ $\PP(Y_t>0,\, \forall\;t\in\RR_+)=1$.
\ Note also that, since \ $X$ \ has continuous trajectories almost surely, by the proof of Lemma \ref{LEMMA_MLE_likelihood},
 we have \ $\PP\Big(\int_0^T\frac{X_s^2}{Y_s}\,\dd s\in[0,\infty)\Big)=1$ \ for all \ $T>0$.
\ Further, for any \ $\omega\in\Omega$, \ $\int_0^T \frac{X_s^2(\omega)}{Y_s(\omega)}\,\dd s =0$ \ holds if and only if
 \ $X_s(\omega) = 0$ \ for almost every \ $s\in[0,T]$.
\ Using that \ $X$ \ has continuous trajectories almost surely, we have
 \[
   \PP\left(\int_0^T \frac{X_s^2}{Y_s}\,\dd s =0 \right)>0
 \]
 holds if and only if \ $\PP(X_s = 0,\; \forall \,s\in[0,T])>0$.
\ Due to \ $a\geq \frac{1}{2}$ \ there does not exist a constant \ $c\in\RR$ \ such that
 \ $\PP(X_s=c,\; \forall\,s\in[0,T])>0$.
\ Indeed, if \ $c\in\RR$ \ is such that \ $\PP(X_s=c,\; \forall\,s\in[0,T])>0$, \ then using \eqref{help2}, we have
 \[
    c = \ee^{-\theta s}\left( c + m\int_0^s \ee^{\theta u}\,\dd u + \int_0^s \ee^{\theta u}\sqrt{Y_u}\,\dd B_u\right),
        \qquad s\in[0,T],
 \]
 on the event \ $\{ \omega\in\Omega : X_s(\omega)=c,\; \forall\,s\in[0,T] \}$.
\ In case \ $\theta\ne 0$, \ the process
 \[
    \left( \int_0^s \ee^{\theta u}\sqrt{Y_u}\,\dd B_u \right)_{s\in[0,T]}
 \]
 would be equal to the deterministic process \ $\left( (c - m/\theta)(\ee^{\theta s} - 1)\right)_{s\in[0,T]}$ \ on the event
 \ $\{ \omega\in\Omega : X_s(\omega)=c,\; \forall\,s\in[0,T] \}$ \ having positive probability.
Since the quadratic variation of the deterministic process \ $\left( (c - m/\theta)(\ee^{\theta s} - 1)\right)_{s\in[0,T]}$ \
 is the identically zero process (due to the fact that the quadratic variation process is a process starting from \ $0$ \ almost surely),
 the quadratic variation of \ $\left( \int_0^s \ee^{\theta u}\sqrt{Y_u}\,\dd B_u \right)_{s\in[0,T]}$ \
 should be identically zero on the event \ $\{ \omega\in\Omega : X_s(\omega)=c,\; \forall\,s\in[0,T]\}$.
\ This would imply that \ $\int_0^s \ee^{2\theta u} Y_u\,\dd u=0$  \ for all \ $s\in[0,T]$
 \ on the event \ $\{ \omega\in\Omega : X_s(\omega)=c,\; \forall\,s\in[0,T]\}$.
\ Using the almost sure continuity and non-negativeness of the sample paths of \ $Y$,  we have
 \ $Y_s=0$  \ for all \ $s\in[0,T]$ \ on the event
 \[
  \{ \omega\in\Omega : X_s(\omega)=c,\; \forall\,s\in[0,T]\}
     \cap \{ \omega\in\Omega : (Y_s(\omega))_{s\in[0,T]}\;\; \text{is continuous}\}
 \]
 having positive probability.
If \ $\theta=0$, \ then replacing \ $\left( (c - m/\theta)(\ee^{\theta s} - 1)\right)_{s\in[0,T]}$ \ by
 \ $(-ms)_{s\in[0,T]}$, \ one can arrive at the same conclusion.
Hence \ $\PP(\inf\{t\in\RR_+ : Y_t=0\} = 0) > 0$ \ which leads us to a contradiction.
This implies \eqref{help40}.

The above argument also shows that we can make a shortcut by arguing in a little bit different way.
Indeed, if \ $C$ \ is a random variable such that \ $\PP(X_s=C,\; \forall\,s\in[0,T])>0$, \ then on the event
 \ $\{ \omega\in\Omega : X_s(\omega)=C(\omega),\; \forall\,s\in[0,T]\}$, \ the quadratic variation of
 \ $X$ \ would be identically zero.
Since, by the SDE \eqref{2dim_affine}, the quadratic variation of \ $X$ \ is the process
 \ $\left(\int_0^t Y_u\,\dd u\right)_{t\geq 0}$, \ it would imply that \ $\int_0^s Y_u\,\dd u =0$ \ for all
 \ $s\in[0,T]$ \ on the event \ $\{ \omega\in\Omega : X_s(\omega)=C(\omega),\; \forall\,s\in[0,T]\}$.
\ Using the almost sure continuity and non-negativeness of the sample paths of \ $Y$, \ we have
 \ $Y_s=0$ \ for all \ $s\in[0,T]$ \ on the event
 \[
  \{ \omega\in\Omega : X_s(\omega)=C(\omega),\; \forall\,s\in[0,T]\}
     \cap \{ \omega\in\Omega : (Y_s(\omega))_{s\in[0,T]}\;\; \text{is continuous}\}
 \]
 having positive probability.
As before, this leads us to a contradiction.
\proofend

The next lemma is about the existence of
 \ $(\widehat{a}_T^{\mathrm{MLE}},\widehat{b}_T^{\mathrm{MLE}},
     \widehat m_T^{\mathrm{MLE}},\widehat\theta_T^{\mathrm{MLE}})$.

\begin{Lem}\label{LEMMA_MLE_exist_theta_m}
If \ $a\geq \frac{1}{2}$, \ $b,m,\theta\in\RR$, \ and \ $\PP(Y_0>0)=1$, \ then
 \begin{gather}
  \label{help41a}
  \PP\left( \int_0^T Y_s\,\dd s \int_0^T\frac{1}{Y_s}\,\dd s - T^2
            \in (0,\infty) \right) = 1
  \qquad \text{for all \ $T>0$,}\\
   \label{help41}
  \PP\left( \int_0^T \frac{X_s^2}{Y_s}\,\dd s \int_0^T\frac{1}{Y_s}\,\dd s
            - \left(\int_0^T \frac{X_s}{Y_s}\,\dd s\right)^2
            \in (0,\infty) \right) = 1
  \qquad \text{for all \ $T>0$,}
 \end{gather}
 and hence there exists a unique MLE
 \ $(\widehat{a}^{\mathrm{MLE}}_T,\widehat{b}^{\mathrm{MLE}}_T,\widehat m^{\mathrm{MLE}}_T,\widehat\theta^{\mathrm{MLE}}_T)$ \
 which has the form given in \eqref{MLE_theta_m_1}, \eqref{MLE_theta_m_2} \eqref{MLE_theta_m_3} and \eqref{MLE_theta_m_4}.
\end{Lem}

\noindent{\bf Proof.}
First note that under the condition \ $a\geq\frac{1}{2}$, \ as it was detailed in the proof of
 Lemma \ref{LEMMA_MLE_exist_theta}, we have \ $\PP(Y_t>0,\, \forall\;t\in\RR_+)=1$.
\ By Cauchy-Schwarz's inequality, we have
 \begin{align}\label{help61_Y}
    \int_0^T Y_s\,\dd s \int_0^T\frac{1}{Y_s}\,\dd s - T^2 \geq 0,\qquad T>0,
 \end{align}
 and hence, using also that \ $Y$ \ has continuous trajectories almost surely and \ $\int_0^T\frac{1}{Y_s}\,\dd s<\infty$ \
 (see the proof of Lemma \ref{LEMMA_MLE_likelihood}), we have
 \[
   \PP\left( \int_0^T Y_s\,\dd s \int_0^T\frac{1}{Y_s}\,\dd s - T^2 \in[0,\infty)\right)=1, \qquad T>0.
 \]
Further, equality holds in \eqref{help61_Y} if and only if \ $K Y_s = L/Y_s$ \ for almost every \ $s\in[0,T]$ \ with some \ $K,L\geq 0$,
 $K^2+L^2>0$ \ ($K$ \ and \ $L$ \ may depend on \ $\omega\in\Omega$ \ and \ $T>0$) \ or equivalently
 \ $K Y_s^2 = L$ \ for almost every \ $s\in[0,T]$ \ with some \ $K,L\geq 0$, $K^2+L^2>0$.
\ Note that if \ $K$ \ were \ $0$, \ then \ $L$ \ would be \ $0$, \ too, hence \ $K$ \ can not be \ $0$ \ implying that
 equality holds in \eqref{help61_Y} if and only if \ $Y_s^2 = L/K$ \ for almost every \ $s\in[0,T]$ \ with some
 \ $K>0$ \ and \ $L\geq 0$ \ ($K$ \ and \ $L$ \ may depend on \ $\omega\in\Omega$ \ and \ $T$).
\ Using that \ $Y$ \ has continuous trajectories almost surely, we have
 \begin{align}\label{help62_Y}
      \PP\left( \int_0^T Y_s\,\dd s \int_0^T\frac{1}{Y_s}\,\dd s - T^2 =0\right)>0
 \end{align}
 holds if and only if \ $\PP(Y_s^2 = L/K,\; \forall \,s\in[0,T])>0$ \ with some random variables \ $K$ \ and \ $L$ \ such that
 \ $\PP(K>0,L\geq 0)=1$ \ ($K$ \ and \ $L$ \ may depend on \ $T$).
\ Similarly, as it was explained at the end of the proof of Lemma \ref{LEMMA_MLE_exist_theta}, this implies that the quadratic variation of
 the process \ $(Y_s^2)_{s\in[0,T]}$ \ would be identically zero on the event
 \ $\{\omega\in\Omega : Y_s^2(\omega) = L(\omega)/K(\omega),\; \forall \,s\in[0,T]\}$ \ having positive probability.
Since, by It\^{o}'s formula,
 \[
   \dd Y_t^2 = 2Y_t\,\dd Y_t + Y_t\,\dd t
             = (2Y_t(a-bY_t) + Y_t)\,\dd t + 2Y_t\sqrt{Y_t}\,\dd L_t,\quad t\geq 0,
 \]
 the quadratic variation of \ $(Y^2_t)_{t\geq 0}$ \ is the process \ $\left(\int_0^t 4Y_u^3\,\dd u\right)_{t\geq 0}$.
\ If \eqref{help62_Y} holds, then \ $\int_0^s 4Y_u^3\,\dd u = 0$ \ for all \ $s\in[0,T]$ \ on the event
 \ $\{\omega\in\Omega : Y_s^2(\omega) = L(\omega)/K(\omega),\; \forall \,s\in[0,T]\}$ \ having positive probability.
Using the almost sure continuity and non-negativeness of the sample paths of \ $Y$, \ we have
 \ $Y_s= 0$ \ for all \ $s\in[0,T]$ \ on the event
 \[
  \{\omega\in\Omega : Y_s^2(\omega) = L(\omega)/K(\omega), \; \forall \,s\in[0,T]\}
     \cap\{ \omega\in\Omega : (Y_t(\omega))_{t\geq 0} \;\; \text{is continuous}\}
 \]
 having positive probability.
This yields us to a contradiction since \ $\PP(Y_t>0, \forall\; t\geq 0)=1$, \ implying \eqref{help41a}.

Now we turn to prove \eqref{help41}.
By Cauchy-Schwarz's inequality, we have
 \begin{align}\label{help61}
    \int_0^T \frac{X_s^2}{Y_s}\,\dd s \int_0^T\frac{1}{Y_s}\,\dd s - \left(\int_0^T \frac{X_s}{Y_s}\,\dd s\right)^2
       \geq 0,\qquad T>0,
 \end{align}
 and hence, since \ $X$ \ has continuous trajectories almost surely, by the proof of Lemma \ref{LEMMA_MLE_likelihood},
 we have
 \[
   \PP\left( \int_0^T \frac{X_s^2}{Y_s}\,\dd s \int_0^T\frac{1}{Y_s}\,\dd s - \left(\int_0^T \frac{X_s}{Y_s}\,\dd s\right)^2
              \in[0,\infty)\right)=1, \qquad T>0.
 \]
Further, equality holds in \eqref{help61} if and only if \ $K X_s^2/Y_s = L/Y_s$ \ for almost every \ $s\in[0,T]$ \ with some \ $K,L\geq 0$,
 $K^2+L^2>0$ \ ($K$ \ and \ $L$ \ may depend on \ $\omega\in\Omega$ \ and \ $T>0$) \ or equivalently
 \ $K X_s^2 = L$ \ for almost every \ $s\in[0,T]$ \ with some \ $K,L\geq 0$, $K^2+L^2>0$.
\ Note that if \ $K$ \ were \ $0$, \ then \ $L$ \ would be \ $0$, \ too, hence \ $K$ \ can not be \ $0$ \ implying that
 equality holds in \eqref{help61} if and only if \ $X_s^2 = L/K$ \ for almost every \ $s\in[0,T]$ \ with some
 \ $K>0$ \ and \ $L\geq 0$ \ ($K$ \ and \ $L$ \ may depend on \ $\omega\in\Omega$ \ and \ $T$).
\ Using that \ $X$ \ has continuous trajectories almost surely, we have
 \begin{align}\label{help62}
      \PP\left( \int_0^T \frac{X_s^2}{Y_s}\,\dd s \int_0^T\frac{1}{Y_s}\,\dd s - \left(\int_0^T \frac{X_s}{Y_s}\,\dd s\right)^2 =0\right)>0
 \end{align}
 holds if and only if \ $\PP(X_s^2 = L/K,\; \forall \,s\in[0,T])>0$ \ with some random variables \ $K$ \ and \ $L$ \ such that
 \ $\PP(K>0,L\geq 0)=1$ \ ($K$ \ and \ $L$ \ may depend on \ $T$).
\ Similarly, as it was explained at the end of the proof of Lemma \ref{LEMMA_MLE_exist_theta}, this implies that the quadratic variation of
 the process \ $(X_s^2)_{s\in[0,T]}$ \ would be identically zero on the event
 \ $\{\omega\in\Omega : X_s^2(\omega) = L(\omega)/K(\omega),\; \forall \,s\in[0,T]\}$ \ having positive probability.
Since, by It\^{o}'s formula,
 \[
   \dd X_t^2 = 2X_t\,\dd X_t + Y_t\,\dd t
             = (2X_t(m-\theta X_t) + Y_t)\,\dd t + 2X_t\sqrt{Y_t}\,\dd B_t,\quad t\geq 0,
 \]
 the quadratic variation of \ $(X^2_t)_{t\geq 0}$ \ is the process \ $\left(\int_0^t 4X_u^2Y_u\,\dd u\right)_{t\geq 0}$.
\ If \eqref{help62} holds, then \ $\int_0^s 4X_u^2Y_u\,\dd u = 0$ \ for all \ $s\in[0,T]$ \ on the event
 \ $\{\omega\in\Omega : X_s^2(\omega) = L(\omega)/K(\omega),\; \forall \,s\in[0,T]\}$ \ having positive probability.
Using the almost sure continuity and non-negativeness of the sample paths of \ $X^2 Y$, \ we have
 \ $X_s^2Y_s = 0$ \ for all \ $s\in[0,T]$ \ on the event
 \[
  \{\omega\in\Omega : X_s^2(\omega) = L(\omega)/K(\omega), \; \forall \,s\in[0,T]\}
     \cap\{ \omega\in\Omega : (X_t^2(\omega)Y_t(\omega))_{t\geq 0} \;\; \text{is continuous}\}
     =: A
 \]
 having positive probability.
Since \ $\PP(Y_t>0,\; \forall\,t\geq 0) = 1$, \ we have \ $X_s = 0$ \ for all \ $s\in[0,T]$ \ on the event \ $A$ \ having
 positive probability.
Repeating the argument given in the proof of Lemma \ref{LEMMA_MLE_exist_theta}, we arrive at a contradiction.
This implies \eqref{help41}.
\proofend

\section{Existence and uniqueness of least squares estimator}
\label{section_EULSE}

Studying LSE for the model \eqref{2dim_affine}, the parameters \ $a>0$ \ and \ $b\in\RR$ \
 will be not supposed to be known.
However, we will not consider the LSEs of \ $a$ \ and \ $b$, \ we will focus only on the LSEs
 of \ $m$ \ and \ $\theta$, \ since we would like use a continuous time observation only for the process \ $X$, \ and
 not for the process \ $(Y,X)$, \ studying LSEs.

First we give a motivation for the LSE based on continuous time observations using the form
 of the LSE based on discrete time low frequency observations.

Let us suppose that \ $m\in\RR$ \ is known \ ($a>0$ \ and \ $b\in\RR$ \ are not supposed to be known).
The LSE of \ $\theta$ \ based on the discrete time observations
 \ $X_i$, $i=0,1,\ldots,n$, \ can be obtained by solving the following extremum problem
 \[
   \widetilde\theta_n^{\textrm{LSE,D}} := \argmin_{\theta\in\RR} \sum_{i=1}^n (X_i - X_{i-1} - (m-\theta X_{i-1}))^2.
 \]
Here in the notation \ $\widetilde\theta_n^{\textrm{LSE,D}}$ \ the letter \ $D$ \ refers to discrete time observations,
 and we note that \ $X_0$ \ denotes an observation for the second coordinate of the initial value of the process \ $(Y,X)$.
\ This definition of LSE of \ $\theta$ \ can be considered as the corresponding one given in Hu and Long \cite[formula (1.2)]{HuLon2}
 for generalized Ornstein-Uhlenbeck processes driven by \ $\alpha$-stable motions, see also Hu and Long
 \cite[formula (3.1)]{HuLon3}.
For a motivation of the LSE of \ $\theta$ \ based on the discrete observations \ $X_i$, \ $i \in \{0, 1, \ldots, n\}$,
 \ see Remark 3.4 in Barczy et al. \cite{BarDorLiPap}.
Further, by Barczy et al. \cite[formula (3.5)]{BarDorLiPap},
 \begin{align}\label{help29}
  \widetilde\theta_n^{\textrm{LSE,D}}
                   = - \frac{\sum_{i=1}^n (X_i - X_{i-1})X_{i-1}  -  m\left(\sum_{i=1}^n X_{i-1}\right)}
                            {\sum_{i=1}^n X_{i-1}^2}
 \end{align}
 provided that \ $\sum_{i=1}^n X_{i-1}^2>0$.
\ Motivated by \eqref{help29}, the LSE of \ $\theta$ \ based on the continuous time observations \ $(X_t)_{t\in[0,T]}$ \
 is defined by
 \begin{align}\label{LSE_theta}
   \widetilde\theta_T^{\mathrm{LSE}} := -\frac{\int_0^T X_s\,\dd X_s - m\int_0^T X_s\,\dd s}
                                  {\int_0^T X_s^2\,\dd s},
 \end{align}
 provided that \ $\int_0^T X_s^2\,\dd s>0$, \ and using the SDE \eqref{2dim_affine} we have
 \begin{align}\label{help27}
   \widetilde\theta_T^{\mathrm{LSE}} - \theta = - \frac{\int_0^T X_s\sqrt{Y_s}\,\dd B_s}{\int_0^T X_s^2\,\dd s},
 \end{align}
 provided that \ $\int_0^T X_s^2\,\dd s>0$.

Let us suppose that the parameters \ $a>0$ \ and \ $b, m\in\RR$ \ are not known.
The LSE of \ $(m,\theta)$ \ based on the discrete time observations \ $X_i$, $i=0,1,\ldots,n$, \ can be obtained
 by solving the following extremum problem
 \[
   (\widehat m_n^{\textrm{LSE,D}},\widehat\theta_n^{\textrm{LSE,D}})
        := \argmin_{(\theta,m)\in\RR^2} \sum_{i=1}^n (X_i - X_{i-1} - (m-\theta X_{i-1}))^2.
 \]
By Barczy et al. \cite[formulas (3.27) and (3.28)]{BarDorLiPap},
  \begin{align}\label{help31}
  \widehat m_n^{\textrm{LSE,D}}
       =  \frac{\sum_{i=1}^n X_{i-1}^2 \sum_{i=1}^n (X_i-X_{i-1})
                   -  \sum_{i=1}^n X_{i-1} \sum_{i=1}^n (X_i-X_{i-1})X_{i-1}}
                            {n\sum_{i=1}^n X_{i-1}^2 - \left(\sum_{i=1}^n X_{i-1}\right)^2}
 \end{align}
 and
 \begin{align}\label{help30}
  \widehat\theta_n^{\textrm{LSE,D}}
        =  \frac{ \sum_{i=1}^n X_{i-1} \sum_{i=1}^n (X_i-X_{i-1}) - n\sum_{i=1}^n (X_i - X_{i-1})X_{i-1} }
                            {n\sum_{i=1}^n X_{i-1}^2 - \left(\sum_{i=1}^n X_{i-1}\right)^2},
 \end{align}
  provided that \ $n\sum_{i=1}^n X_{i-1}^2 - \left(\sum_{i=1}^n X_{i-1}\right)^2>0$.
\ Motivated by \eqref{help31} and \eqref{help30}, the LSE of \ $(m,\theta)$ \ based on the continuous time observations
 \ $(X_t)_{t\in[0,T]}$ \ is defined by
 \begin{align}\label{LSE_theta_m_2}
  \widehat m_T^{\mathrm{LSE}}
     & := \frac{(X_T-X_0)\int_0^T X_s^2\,\dd s - \left(\int_0^T X_s\,\dd s\right)\left(\int_0^T X_s\,\dd X_s\right)}
            {T\int_0^T X_s^2\,\dd s - \left(\int_0^T X_s\,\dd s\right)^2},\\ \label{LSE_theta_m_1}
   \widehat\theta_T^{\mathrm{LSE}}
     & := \frac{(X_T-X_0)\int_0^T X_s\,\dd s - T \int_0^T X_s\,\dd X_s}
            {T\int_0^T X_s^2\,\dd s - \left(\int_0^T X_s\,\dd s\right)^2},
 \end{align}
 provided that \ $T\int_0^T X_s^2\,\dd s - \left(\int_0^T X_s\,\dd s\right)^2>0$.
\ Note that, by Cauchy-Schwarz's inequality, \ $T\int_0^T X_s^2\,\dd s - \left(\int_0^T X_s\,\dd s\right)^2\geq 0$, \
 and \ $T\int_0^T X_s^2\,\dd s - \left(\int_0^T X_s\,\dd s\right)^2>0$ \ yields that \ $\int_0^T X_s^2\,\dd s >0$.
\ Then
 \[
  \begin{bmatrix}
    \widehat m_T^{\mathrm{LSE}} \\
    \widehat\theta_T^{\mathrm{LSE}} \\
  \end{bmatrix}
   = \begin{bmatrix}
       T & -\int_0^T X_s\,\dd s \\
        -\int_0^T X_s\,\dd s & \int_0^T X_s^2\,\dd s \\
     \end{bmatrix}^{-1}
     \begin{bmatrix}
        X_T - X_0 \\
        -\int_0^T X_s\,\dd X_s\\
     \end{bmatrix},
 \]
 and using the SDE \eqref{2dim_affine} one can check that
 \begin{align*}
  \begin{bmatrix}
    \widehat m_T^{\mathrm{LSE}} - m \\
    \widehat\theta_T^{\mathrm{LSE}} - \theta \\
  \end{bmatrix}
  =
   \begin{bmatrix}
     T & -\int_0^T X_s\,\dd s \\
     -\int_0^T X_s\,\dd s & \int_0^T X_s^2\,\dd s \\
   \end{bmatrix}^{-1}
   \begin{bmatrix}
     \int_0^T \sqrt{Y_s}\,\dd B_s  \\
      - \int_0^T X_s\sqrt{Y_s}\,\dd B_s\\
   \end{bmatrix},
 \end{align*}
 and hence
 \begin{align}\label{help45}
   \widehat m_T^{\mathrm{LSE}} - m
     =  \frac{-\int_0^T X_s\,\dd s \int_0^T X_s\sqrt{Y_s}\,\dd B_s + \int_0^T X_s^2\,\dd s  \int_0^T \sqrt{Y_s}\,\dd B_s }
             {T\int_0^T X_s^2\,\dd s - \left(\int_0^T X_s\,\dd s\right)^2},
 \end{align}
 and
 \begin{align}\label{help44}
   \widehat\theta_T^{\mathrm{LSE}} - \theta
     =  \frac{-T\int_0^T X_s\sqrt{Y_s}\,\dd B_s + \int_0^T X_s\,\dd s  \int_0^T \sqrt{Y_s}\,\dd B_s }
             {T\int_0^T X_s^2\,\dd s - \left(\int_0^T X_s\,\dd s\right)^2},
 \end{align}
 provided that \ $T\int_0^T X_s^2\,\dd s - \left(\int_0^T X_s\,\dd s\right)^2>0$.

\begin{Rem}
For the stochastic integral \ $\int_0^T X_s\,\dd X_s$ \ in \eqref{LSE_theta},
 \eqref{LSE_theta_m_2} and \eqref{LSE_theta_m_1}, we have
 \[
   \sum_{i=1}^{\lfloor nT\rfloor } X_{\frac{i-1}{n}}
             \big(X_{\frac{i}{n}} - X_{\frac{i-1}{n}} \big)
   \stoch \int_0^T X_s\,\dd X_s \qquad \text{as \ $n\to\infty$,}
 \]
 following from Proposition I.4.44 in Jacod and Shiryaev \cite{JSh}.
For more details, see Remark \ref{Riemann}.
\proofend
\end{Rem}

The next lemma is about the existence of \ $\widetilde \theta_T^{\mathrm{LSE}}$ \ (supposing that \ $m\in\RR$ \ is known,
 but \ $a>0$ \ and \ $b$ \ are unknown).

\begin{Lem}\label{LEMMA_LSE_exist_theta}
If \ $a>0$, \ $b,m,\theta\in\RR$, \ and \ $\PP(Y_0>0)=1$, \ then
 \begin{align}\label{help39}
  \PP\left( \int_0^T X_s^2\,\dd s \in (0,\infty)\right)=1\quad \text{for all \ $T>0$,}
 \end{align}
 and hence there exists a unique LSE \ $\widetilde\theta^{\mathrm{LSE}}_T$ \ which has the form given in \eqref{LSE_theta}.
\end{Lem}

\noindent{\bf Proof.}
First note that under the condition on the parameters, for all \ $y_0>0$, \ we have
 \[
   \PP(\inf\{ t\in\RR_+ : Y_t = 0\} >0 \mid Y_0 = y_0)=1,
 \]
 and hence, by the law of total probability, \ $\PP(\inf\{ t\in\RR_+ : Y_t = 0\} >0 )=1$.
\ Since \ $X$ \ has continuous trajectories almost surely, we readily have
 \[
  \PP\left( \int_0^T X_s^2\,\dd s \in[0,\infty)\right)=1,
   \qquad T>0.
 \]
Observe also that for any \ $\omega\in\Omega$, \
 $\int_0^T X_s^2(\omega)\,\dd s =0$ \ holds if and only if \ $X_s(\omega) = 0$ \ for almost every \ $s\in[0,T]$.
\ Using that \ $X$ \ has continuous trajectories almost surely, we have
 \[
   \PP\left(\int_0^T X_s^2\,\dd s =0 \right)>0
 \]
 holds if and only if \ $\PP(X_s = 0,\; \forall \,s\in[0,T])>0$.
By the end of the proof of Lemma \ref{LEMMA_MLE_exist_theta},
 if \ $\PP(X_s = 0,\, \forall\, s\in[0,T])>0$, \ then \ $Y_s=0$ \ for all \ $s\in[0,T]$ \ on the event
 \[
  \{\omega\in\Omega : X_s(\omega) = 0,\; \forall\, s\in[0,T] \}
    \cap \big\{\omega\in\Omega : (Y_s(\omega))_{s\in[0,T]}  \text{\;\; is continuous}  \big\}
 \]
 having positive probability.
This would yield that \ $\PP(\inf\{t\in\RR_+ : Y_t = 0\} = 0)>0$ \ leading us to a contradiction,
 which implies \eqref{help39}.
\proofend

The next lemma is about the existence of \ $(\widehat m_T^{\mathrm{LSE}},\widehat\theta_T^{\mathrm{LSE}})$.

\begin{Lem}\label{LEMMA_LSE_exist_theta_m}
If \ $a>0$, \ $b,m,\theta\in\RR$, \ and \ $\PP(Y_0>0)=1$, \ then
 \begin{align}\label{help38}
      \PP\left( T\int_0^T X_s^2\,\dd s - \left(\int_0^T X_s\,\dd s\right)^2 \in (0,\infty)\right)=1
         \quad \text{for all \ $T>0$,}
 \end{align}
 and hence there exists a unique LSE \ $(\widehat m^{\mathrm{LSE}}_T,\widehat\theta^{\mathrm{LSE}}_T)$ \
 which has the form given in \eqref{LSE_theta_m_2} and \eqref{LSE_theta_m_1}.
\end{Lem}

\noindent{\bf Proof.}
Just as in the proof of Lemma \ref{LEMMA_LSE_exist_theta}, we have \ $\PP(\inf\{t\in\RR_+ : Y_t = 0\} > 0)=1$.
\ By Cauchy-Schwarz's inequality, we have
 \begin{align}\label{help36}
    T\int_0^T X^2_s\,\dd s - \left(\int_0^T X_s\,\dd s\right)^2 \geq 0,\qquad T>0,
 \end{align}
 and hence, since \ $X$ \ has continuous trajectories almost surely, we readily have
 \[
    \PP\left( T\int_0^T X^2_s\,\dd s - \left(\int_0^T X_s\,\dd s\right)^2 \in[0,\infty) \right)=1,
    \qquad T>0.
 \]
Further, equality holds in \eqref{help36} if and only if \ $K X_s^2 = L$ \ for almost every \ $s\in[0,T]$ \
 with some \ $K,L\geq0$, $K^2+L^2>0$ \ ($K$ \ and \ $L$ \ may depend on \ $\omega\in\Omega$ \ and \ $T$).
\ Note that if \ $K$ \ were \ $0$, \ then \ $L$ \ would be \ $0$, \ too, hence \ $K$ \ can not be \ $0$ \ implying that
 equality holds in \eqref{help36} if and only if \ $X_s^2 = L/K$ \ for almost every \ $s\in[0,T]$ \ with some
 \ $K>0$ \ and \ $L\geq 0$ \ ($K$ \ and \ $L$ \ may depend on \ $\omega\in\Omega$ \ and \ $T$).
\ Using that \ $X$ \ has continuous trajectories almost surely, we have
 \begin{align}\label{help76}
      \PP\left( T\int_0^T X_s^2\,\dd s - \left(\int_0^T X_s\,\dd s\right)^2 =0\right)>0
 \end{align}
 holds if and only if \ $\PP(X_s^2 = L/K,\; \forall \,s\in[0,T])>0$ \ with some random variables \ $K$ \ and \ $L$ \ such that
 \ $\PP(K>0,L\geq 0)=1$ \ ($K$ \ and \ $L$ \ may depend on \ $T$).
\ Using again that \ $X$ \ has continuous trajectories almost surely, we have \eqref{help76} holds if and only if
 \ $\PP(X_s = C,\; \forall \,s\in[0,T])>0$ \ with some random variable \ $C$ \ (note that \ $C=\sqrt{L/K}$ \ if
 \ $X$ \ is non-negative and \ $C=-\sqrt{L/K}$ \ if \ $X$ \ is negative).
This leads us to a contradiction by the end of the proof of Lemma \ref{LEMMA_MLE_exist_theta}.
Indeed, if \ $\PP(X_s = C,\;\forall\, s\in[0,T])>0$ \ with some random variable \ $C$, \ then, by the end of the proof
 of Lemma \ref{LEMMA_MLE_exist_theta}, we have \ $Y_s = 0$ \ for all \ $s\in[0,T]$ \ on the event
 \[
  \Big\{\omega\in\Omega : X_s(\omega) = C(\omega),\; \forall\, s\in[0,T] \Big\}
    \cap \big\{\omega\in\Omega : (Y_s(\omega))_{s\in[0,T]}  \text{\;\; is continuous}  \big\}
 \]
 having positive probability.
This would yield that \ $\PP(\inf\{t\in\RR_+ : Y_t = 0\} = 0)>0$ \ leading us to a contradiction,
 which implies \eqref{help38}.
\proofend

\section{Consistency of maximum likelihood estimator}
\label{section_CMLE}

\begin{Thm}\label{Thm_MLE_cons_theta}
If \ $a\geq \frac{1}{2}$, \ $b>0$, \ $m\in\RR$, \ $\theta>0$,
 \ and \ $\PP(Y_0>0)=1$, \ then the MLE of \ $\theta$ \ is
 strongly consistent: \ $\PP\left(\lim_{T\to\infty} {\widetilde\theta}^{\mathrm{MLE}}_T = \theta\right)=1$.
\end{Thm}

\noindent{\bf Proof.}
By Lemma \ref{LEMMA_MLE_exist_theta}, there exists a unique MLE \ ${\widetilde\theta}^{\mathrm{MLE}}_T $ \ of \ $\theta$ \
 which has the form given in \eqref{MLE_theta}.
Further, by \eqref{help52},
 \[
  {\widetilde\theta}^{\mathrm{MLE}}_T  - \theta
  = - \frac{\int_0^T \frac{X_s}{\sqrt{Y_s}}\,\dd B_s}
           {\int_0^T \frac{X_s^2}{Y_s}\,\dd s},
  \qquad T > 0 .
 \]
The strong consistency of \ ${\widetilde\theta}^{\mathrm{MLE}}_T$ \ will follow from
 a strong law of large numbers for continuous local martingales (see, e.g., Theorem \ref{DDS_stoch_int})
 provided that we check that
  \begin{align} \label{help82}
   & \PP\left(\int_0^T \frac{X_s^2}{Y_s}\,\dd s < +\infty \right)=1,\qquad T>0,\\
   &\PP\left(\lim_{T\to\infty} \int_0^T \frac{X_s^2}{Y_s}\,\dd s = +\infty \right)=1.\label{help10}
 \end{align}
Since \ $(Y,X)$ \ has continuous trajectories almost surely, we have \eqref{help82}.
Next we check \eqref{help10}.
By Theorem \ref{Thm_ergodic12}, we have
 \[
    \lim_{T\to\infty} \frac{1}{T} \int_0^T \frac{X_s^2}{1+Y_s}\,\dd s = \EE\left(\frac{X_\infty^2}{1+Y_\infty}\right)
    \qquad \text{a.s.},
 \]
 where
 \[
    \EE\left(\frac{X_\infty^2}{1+Y_\infty}\right) \leq \EE(X_\infty^2)
              = \frac{a\theta + 2bm^2}{2b\theta^2} <\infty.
 \]
Next we check that \ $\EE(X_\infty^2/(1+Y_\infty))>0$. \
Since \ $X_\infty^2/(1+Y_\infty)\geq 0$, \ we have \ $\EE(X_\infty^2/(1+Y_\infty))=0$ \ holds if and only if
 \ $\PP(X_\infty^2/(1+Y_\infty) = 0)=1$ \ or equivalently \ $\PP(X_\infty = 0)=1$ \ which leads us to a contradiction
 since \ $X_\infty$ \ is absolutely continuous by Theorem \ref{Thm_ergodic12}.
Hence we have
 \[
    \PP\left(  \lim_{T\to\infty} \int_0^T \frac{X_s^2}{1+Y_s}\,\dd s  = \infty \right)=1,
 \]
which yields \eqref{help10}.
\proofend

\begin{Thm}\label{Thm_MLE_cons_theta_m}
If \ $a>\frac{1}{2}$, \ $b>0$, \ $m\in\RR$, \ $\theta>0$, \ and \ $\PP(Y_0>0)=1$,
 \ then the MLE of \ $(a,b,m,\theta)$ \ is strongly consistent:
 \ $\PP\left(\lim_{T\to\infty}
              ({\widehat a}^{\mathrm{MLE}}_T,{\widehat b}^{\mathrm{MLE}}_T,
               {\widehat m}^{\mathrm{MLE}}_T,{\widehat\theta}^{\mathrm{MLE}}_T)
             = (a,b,m,\theta)\right)
    =1$.
\end{Thm}

\noindent{\bf Proof.}
By Lemma \ref{LEMMA_MLE_exist_theta_m}, there exists a unique MLE
 \ $({\widehat a}^{\mathrm{MLE}}_T,{\widehat b}^{\mathrm{MLE}}_T,
          {\widehat m}^{\mathrm{MLE}}_T,{\widehat\theta}^{\mathrm{MLE}}_T)$ \
 of \ $(a,b,m,\theta)$ \ which has the form given in
 \eqref{MLE_theta_m_1}, \eqref{MLE_theta_m_2}, \eqref{MLE_theta_m_3} and
 \eqref{MLE_theta_m_4}.
First we check that \ $\PP(\lim_{T\to\infty} {\widehat\theta}^{\mathrm{MLE}}_T = \theta)=1$.
\ By \eqref{help48}, using also that
 \ $\int_0^T \frac{X_s^2}{Y_s}\,\dd s \int_0^T\frac{1}{Y_s}\,\dd s - \left(\int_0^T \frac{X_s}{Y_s}\,\dd s\right)^2 >0$ \
 yields \ $\int_0^T \frac{X_s^2}{Y_s}\,\dd s \int_0^T\frac{1}{Y_s}\,\dd s>0$, \ we have
 \begin{align}\label{help57}
   {\widehat\theta}^{\mathrm{MLE}}_T - \theta
     = \frac{\frac{\frac{1}{T}\int_0^T \frac{X_s}{Y_s}\,\dd s}{\frac{1}{T}\int_0^T \frac{X_s^2}{Y_s}\,\dd s}
             \cdot \frac{ \int_0^T \frac{1}{\sqrt{Y_s}}\,\dd B_s}{\int_0^T \frac{1}{Y_s}\,\dd s}
             - \frac{ \int_0^T \frac{X_s}{\sqrt{Y_s}}\,\dd B_s }{\int_0^T \frac{X_s^2}{Y_s}\,\dd s}}
         {1 - \frac{\left(  \frac{1}{T}\int_0^T \frac{X_s}{Y_s}\,\dd s \right)^2}{\frac{1}{T}\int_0^T \frac{X_s^2}{Y_s}\,\dd s
                            \cdot \frac{1}{T}\int_0^T \frac{1}{Y_s}\,\dd s }}
  \qquad \text{a.s.}
 \end{align}
 due to \eqref{help41}.
Next, we show that \ $\EE(1/ Y_\infty)<\infty$, \ $\EE(X_\infty/ Y_\infty)<\infty$ \ and
 \ $\EE(X_\infty^2/ Y_\infty)<\infty$, \ which, by part (iii) of Theorem \ref{Thm_ergodic12}, will imply that
 \begin{align}\label{help55}
   &\PP\left( \lim_{T\to\infty} \frac{1}{T}\int_0^T \frac{1}{Y_s}\,\dd s = \EE\left(\frac{1}{Y_\infty}\right)\right)=1,\\ \label{help56}
   &\PP\left( \lim_{T\to\infty} \frac{1}{T}\int_0^T \frac{X_s}{Y_s}\,\dd s = \EE\left(\frac{X_\infty}{Y_\infty}\right)\right)=1,\\
   & \label{help42}
    \PP\left( \lim_{T\to\infty} \frac{1}{T}\int_0^T \frac{X_s^2}{Y_s}\,\dd s = \EE\left(\frac{X_\infty^2}{Y_\infty}\right)\right)=1.
 \end{align}
We only prove \ $\EE(X_\infty^2/ Y_\infty)<\infty$ \ noting that \ $\EE(1/ Y_\infty)<\infty$ \ and
 \ $\EE(X_\infty/ Y_\infty)<\infty$ \ can be checked in the very same way.
First note that, by Theorem \ref{Thm_ergodic12}, \ $\PP(Y_\infty >0)=1$, \ and hence the random variable
 \ $X_\infty^2/ Y_\infty$ \ is well-defined with probability one.
By H\"older's inequality, for all \ $p,q>0$ \ with \ $\frac{1}{p}+\frac{1}{q}=1$, \ we have
 \[
  \EE\left(\frac{X_\infty^2}{Y_\infty}\right)
     \leq (\EE(X_\infty^{2p}))^{1/p}  \left( \EE\left(\frac{1}{Y_\infty^q}\right)\right)^{1/q}.
 \]
Since, by Theorem \ref{Thm_ergodic12}, \ $\EE(\vert X_\infty\vert^n)<\infty$ \ for all \ $n\in\NN$,
 \ to prove that \ $\EE\left(\frac{X_\infty^2}{Y_\infty}\right)<\infty$ \ it is enough to check that
 there exists some \ $\varepsilon>0$ \ such that \ $\EE\left(\frac{1}{Y_\infty^{1+\varepsilon}}\right)<\infty$.
\ Using that, by Theorem \ref{Thm_ergodic12}, \ $Y_\infty$ \ has Gamma distribution with density function
 \ $\frac{(2b)^{2a}}{\Gamma(2a)}x^{2a-1}\ee^{-2bx}\mathbf 1_{\{x>0\}}$, \ we have, for any \ $\varepsilon>0$,
  \begin{align*}
    \EE\left(\frac{1}{Y_\infty^{1+\varepsilon}}\right)
       = \int_0^\infty \frac{1}{x^{1+\varepsilon}}
                       \cdot\frac{(2b)^{2a}}{\Gamma(2a)} x^{2a-1}\ee^{-2bx}\,\dd x
       = \frac{(2b)^{2a}}{\Gamma(2a)}
         \int_0^\infty  x^{2a-2-\varepsilon}\ee^{-2bx}\,\dd x.
  \end{align*}
Due to our assumption \ $a>1/2$, \ one can choose an \ $\varepsilon$ \ such that \ $\max(0,2a-2)<\varepsilon<2a-1$,
 \ and hence for all \ $M>0$,
 \begin{align*}
  \int_0^\infty  x^{2a-2-\varepsilon}\ee^{-2bx}\,\dd x
     & \leq \int_0^M x^{2a-2-\varepsilon} \,\dd x
           + M^{2a-2-\varepsilon} \int_M^\infty \ee^{-2bx}\,\dd x \\
     & = \frac{M^{2a-1-\varepsilon}}{2a-1-\varepsilon}
        + M^{2a-2-\varepsilon}
          \lim_{L\to\infty}\frac{\ee^{-2bL} - \ee^{-2bM} }{-2b}\\
     & = \frac{M^{2a-1-\varepsilon}}{2a-1-\varepsilon}
        + M^{2a-2-\varepsilon} \frac{\ee^{-2bM}}{2b}
      <\infty,
 \end{align*}
 which yields that \ $\EE(1/Y_\infty^{1+\varepsilon})<\infty$ \ for \ $\varepsilon$ \ satisfying
 \ $\max(0,2a-2)<\varepsilon<2a-1$.

Further, since \ $\EE(1/ Y_\infty)>0$ \ and \ $\EE(X_\infty^2/ Y_\infty)>0$ \ (due to the absolutely continuity of
  \ $X_\infty$, \ as it was explained in the proof of Theorem \ref{Thm_MLE_cons_theta}),
 \eqref{help55} and \eqref{help42} yield that
 \begin{align*}
    \PP\left( \lim_{T\to\infty} \int_0^T \frac{1}{Y_s}\,\dd s = \infty \right)
       = \PP\left( \lim_{T\to\infty} \int_0^T \frac{X_s^2}{Y_s}\,\dd s = \infty\right)=1.
 \end{align*}
Using \eqref{help57}, \eqref{help55}, \eqref{help56}, \eqref{help42},
 Slutsky's lemma and a strong law of large numbers for continuous local martingales (see, e.g., Theorem \ref{DDS_stoch_int}), we get
 \begin{align*}
  \lim_{T\to\infty} \left({\widehat\theta}^{\mathrm{MLE}}_T - \theta \right)
     = \frac{ \frac{\EE\left(\frac{X_\infty}{Y_\infty}\right)}{\EE\left(\frac{X_\infty^2}{Y_\infty}\right)} \cdot 0 - 0 }
            {1 - \frac{\left(\EE\left(\frac{X_\infty}{Y_\infty}\right)\right)^2}
                      {\EE\left(\frac{X_\infty^2}{Y_\infty}\right)\EE\left(\frac{1}{Y_\infty}\right)} }
    = 0\qquad \text{a.s.,}
 \end{align*}
 where we also used that the denominator is strictly positive.
Indeed, by Cauchy-Schwarz's inequality,
 \begin{align}\label{help59}
  \left(\EE\left(\frac{X_\infty}{Y_\infty}\right)\right)^2
     \leq  \EE\left(\frac{X_\infty^2}{Y_\infty}\right) \EE\left(\frac{1}{Y_\infty}\right),
 \end{align}
 and equality would hold if and only if \ $\PP(KX_\infty^2/Y_\infty = L/Y_\infty)=1$ \ with some \ $K,L\geq 0$, \ $K^2+L^2>0$.
\ By Theorem \ref{Thm_ergodic12}, \ $(Y_\infty,X_\infty)$ \ is absolutely continuous and then
 \ $\PP(X_\infty = c)=0$ \ for all \ $c\in\RR$.
\ This implies that \eqref{help59} holds with strict inequality.

Similarly, by \eqref{help47a}, \eqref{help47b} and \eqref{help47}, we have
 \begin{gather*}
  {\widehat{a}}^{\mathrm{MLE}}_T - a
  =\frac{\frac{\int_0^T \frac{1}{\sqrt{Y_s}}\,\dd L_s}
              {\int_0^T \frac{1}{Y_s}\,\dd s}
         -\frac{1}{\frac{1}{T}\int_0^T \frac{1}{Y_s}\,\dd s}
          \cdot \frac{\int_0^T \sqrt{Y_s}\,\dd L_s}
                     {\int_0^T Y_s\,\dd s}}
        {1-\frac{1}
                {\frac{1}{T}\int_0^T Y_s\,\dd s
                 \cdot \frac{1}{T}\int_0^T \frac{1}{Y_s}\,\dd s }}
  \qquad \text{a.s.,} \\
  {\widehat{b}}^{\mathrm{MLE}}_T - b
  =\frac{\frac{1}{\frac{1}{T}\int_0^T Y_s\,\dd s}
         \cdot \frac{\int_0^T\frac{1}{\sqrt{Y_s}}\,\dd L_s}
                    { \int_0^T \frac{1}{Y_s}\,\dd s}
         -\frac{\int_0^T \sqrt{Y_s}\,\dd L_s}{\int_0^T Y_s\,\dd s}}
        {1-\frac{1}
                {\frac{1}{T}\int_0^T Y_s\,\dd s
                 \cdot \frac{1}{T}\int_0^T \frac{1}{Y_s}\,\dd s }}
  \qquad \text{a.s.,} \\
  {\widehat m}^{\mathrm{MLE}}_T - m
  =\frac{\frac{\int_0^T \frac{1}{\sqrt{Y_s}}\,\dd B_s}
              {\int_0^T \frac{1}{Y_s}\,\dd s}
         -\frac{\frac{1}{T}\int_0^T \frac{X_s}{Y_s}\,\dd s}
          {\frac{1}{T}\int_0^T \frac{1}{Y_s}\,\dd s}
          \cdot \frac{\int_0^T\frac{X_s}{\sqrt{Y_s}}\,\dd B_s}
                     {\int_0^T \frac{X_s^2}{Y_s}\,\dd s}}
        {1-\frac{\left(  \frac{1}{T}\int_0^T \frac{X_s}{Y_s}\,\dd s \right)^2}
                {\frac{1}{T}\int_0^T \frac{X_s^2}{Y_s}\,\dd s
                 \cdot \frac{1}{T}\int_0^T \frac{1}{Y_s}\,\dd s }}
  \qquad \text{a.s.,}
 \end{gather*}
 due to \eqref{help41a} and \eqref{help41}.
Since \ $\EE(Y_\infty)<\infty$, \ by part (iii) of Theorem \ref{Thm_ergodic12},
 \begin{align}\label{help79}
  \PP\left( \lim_{T\to\infty} \frac{1}{T} \int_0^T Y_s\,\dd s = \EE(Y_\infty) \right)
    = \PP\left( \lim_{T\to\infty} \int_0^T Y_s\,\dd s = \infty \right) = 1,
 \end{align}
 and then, similarly as before, one can argue that
 \begin{gather*}
  \lim_{T\to\infty} ({\widehat{a}}^{\mathrm{MLE}}_T - a )
  =\frac{0-\frac{1}{\EE\left(\frac{1}{Y_\infty}\right)} \cdot 0}
        {1-\frac{1}
                {\EE\left(Y_\infty\right)
                 \EE\left(\frac{1}{Y_\infty}\right)}}
    = 0\qquad \text{a.s.},\\
  \lim_{T\to\infty} ({\widehat{b}}^{\mathrm{MLE}}_T - b )
  =\frac{\frac{1}{\EE\left( Y_\infty \right)} \cdot 0 - 0}
        {1-\frac{1}
                {\EE\left(Y_\infty\right)
                 \EE\left(\frac{1}{Y_\infty}\right)}}
    = 0\qquad \text{a.s.},\\
  \lim_{T\to\infty} ({\widehat m}^{\mathrm{MLE}}_T - m )
  =\frac{0-\frac{\EE\left(\frac{X_\infty}{Y_\infty}\right)}
                {\EE\left(\frac{1}{Y_\infty}\right)} \cdot 0}
        {1-\frac{\left(\EE\left(\frac{X_\infty}{Y_\infty}\right)\right)^2}
                {\EE\left(\frac{X_\infty^2}{Y_\infty}\right)
                 \EE\left(\frac{1}{Y_\infty}\right)}}
    = 0\qquad \text{a.s.,}
 \end{gather*}
 where we also used that \ $\EE(Y_\infty)\EE\left(\frac{1}{Y_\infty}\right)>1$ \ (which can be checked
  using Cauchy-Schwarz's inequality and the absolute continuity of \ $Y_\infty$).
Using that the intersection of four events with probability one is an event with probability one, we have the assertion.
\proofend

\begin{Rem}
If \ $a=\frac{1}{2}$, \ $b>0$, \ $\theta>0$, \ $m\in\RR$, \ and \ $\PP(Y_0>0)=1$,
 \ then one should find another approach for studying the consistency behaviour of the MLE of \ $(a,b,m,\theta)$, \ since in this case
 \[
    \EE\left(\frac{1}{Y_\infty}\right)
       = \int_0^\infty \frac{2b\ee^{-2bx}}{x}\,\dd x
       =\infty,
 \]
 and hence one cannot use part (iii) of Theorem \ref{Thm_ergodic12}.
In this paper we renounce to consider it.
\proofend
\end{Rem}

\section{Consistency of least squares estimator}
\label{section_CLSE}

\begin{Thm}\label{Thm_LSE_cons_theta}
If \ $a>0$, \ $b>0$, \ $m\in\RR$, $\theta>0$, \ and \ $\PP(Y_0>0)=1$,
 \ then the LSE of \ $\theta$ \ is strongly consistent:
 \ $\PP\left(\lim_{T\to\infty} {\widetilde\theta}^{\mathrm{LSE}}_T = \theta\right)=1$.
\end{Thm}

\noindent{\bf Proof.}
By Lemma \ref{LEMMA_LSE_exist_theta}, there exists a unique \ ${\widetilde\theta}^{\mathrm{LSE}}_T$ \ of \ $\theta$ \
 which has the form given in \eqref{LSE_theta}.
By \eqref{help27}, we have
 \begin{align}\label{help9}
   \widetilde\theta_T^{\mathrm{LSE}} - \theta
     = - \frac{\int_0^T X_s\sqrt{Y_s}\,\dd B_s}{\int_0^T X_s^2\,\dd s}
     = - \frac{\int_0^T X_s\sqrt{Y_s}\,\dd B_s}{\int_0^T X_s^2Y_s\,\dd s}
        \cdot \frac{\frac{1}{T}\int_0^T X_s^2 Y_s\,\dd s}{\frac{1}{T}\int_0^T X_s^2\,\dd s}.
 \end{align}
By Theorem \ref{Thm_ergodic12}, we have
 \begin{align}\label{help80}
   \PP\left(\lim_{T\to\infty} \frac{1}{T}\int_0^T X_s^2 Y_s\,\dd s = \EE\big(X_\infty^2Y_\infty\big) \right)=1
    \quad \text{and} \quad
   \PP\left(\lim_{T\to\infty} \frac{1}{T}\int_0^T X_s^2 \,\dd s = \EE\big(X_\infty^2\big) \right)=1.
 \end{align}
We note that \ $\EE(X_\infty^2Y_\infty) $ \ and \ $\EE(X_\infty^2)$ \ are calculated explicitly in
 Theorem \ref{Thm_ergodic12}.
Note also that \ $\EE(X_\infty^2Y_\infty)$ \ is positive (due to that \ $X_\infty^2Y_\infty$ \ is non-negative and absolutely continuous),
 \ and hence we also have
 \begin{align*}
   \PP\left(\lim_{T\to\infty} \int_0^T X_s^2 Y_s\,\dd s = \infty \right)=1.
 \end{align*}
Then, by a strong law of large numbers for continuous local martingales (see, e.g., Theorem \ref{DDS_stoch_int}), we get
 \[
    \PP\left(\lim_{T\to\infty}\frac{\int_0^T X_s\sqrt{Y_s}\,\dd B_s}{\int_0^T X_s^2Y_s\,\dd s} =0\right)=1,
 \]
 and hence \eqref{help9} yields the assertion.
\proofend

\begin{Thm}\label{Thm_LSE_cons_theta_m}
If \ $a>0$, \ $b>0$, \ $m\in\RR$, \ $\theta>0$, \ and \ $\PP(Y_0>0)=1$,
 \ then the LSE of \ $(m,\theta)$ \ is strongly consistent:
 \ $\PP\left(\lim_{T\to\infty} ({\widehat m}^{\mathrm{LSE}}_T,{\widehat\theta}^{\mathrm{LSE}}_T)= (m,\theta)\right)=1$.
\end{Thm}

\noindent{\bf Proof.}
By Lemma \ref{LEMMA_LSE_exist_theta_m}, there exists a unique LSE \ $(\widehat m^{\mathrm{LSE}}_T,\widehat\theta^{\mathrm{LSE}}_T)$ \
 of \ $(m,\theta)$ \ which has the form given in \eqref{LSE_theta_m_2} and \eqref{LSE_theta_m_1}.
By \eqref{help44}, we have
 \begin{align*}
   {\widehat \theta}^{\mathrm{LSE}}_T - \theta
     = \frac{  -\frac{\frac{1}{T}\int_0^T X_s^2Y_s\,\dd s}{\frac{1}{T}\int_0^T X_s^2 \,\dd s}
             \cdot \frac{ \int_0^T X_s\sqrt{Y_s}\,\dd B_s}{\int_0^T X_s^2Y_s\,\dd s}
             + \frac{1}{T}\int_0^T X_s\,\dd s
               \cdot\frac{\frac{1}{T}\int_0^T Y_s \,\dd s}{\frac{1}{T}\int_0^T  X_s^2 \,\dd s }
               \cdot\frac{ \int_0^T \sqrt{Y_s}\,\dd B_s }{\int_0^T Y_s\,\dd s}}
         {1 - \frac{\left(  \frac{1}{T}\int_0^T X_s \,\dd s \right)^2}{\frac{1}{T}\int_0^T  X_s^2 \,\dd s }}
  \qquad \text{a.s.}
 \end{align*}
 due to \eqref{help38}.
Using \eqref{help79}, \eqref{help80} and that
 \begin{align}\label{help81}
   \PP\left( \lim_{T\to\infty} \frac{1}{T}\int_0^T X_s\,\dd s = \EE(X_\infty) \right)= 1,
 \end{align}
 similarly to the proof of Theorem \ref{Thm_MLE_cons_theta_m}, we get
 \begin{align*}
  \lim_{T\to\infty} ({\widehat\theta}^{\mathrm{LSE}}_T - \theta )
     = \frac{ -\frac{\EE\left( X_\infty^2 Y_\infty\right)}{\EE( X_\infty^2)} \cdot 0
             + \EE(X_\infty)\cdot \frac{\EE(Y_\infty)}{\EE(X_\infty^2)}\cdot 0 }
            {1 - \frac{(\EE( X_\infty))^2}{\EE(X_\infty^2)} }
    = 0\qquad \text{a.s.,}
 \end{align*}
 where for the last step we also used that \ $(\EE( X_\infty))^2 < \EE(X_\infty^2)$ \
 (which holds since there does not exist a constant \ $c\in\RR$ \ such that \ $\PP(X_\infty = c)=1$
  \ due to the fact that \ $X_\infty$ \ is absolutely continuous).

Similarly, by \eqref{help45},
 \begin{align*}
  \lim_{T\to\infty} ({\widehat m}^{\mathrm{LSE}}_T - m )
   &= \lim_{T\to\infty}
       \frac{ - \frac{1}{T}\int_0^T X_s \,\dd s \cdot \frac{\frac{1}{T}\int_0^T X_s^2Y_s\,\dd s}{\frac{1}{T}\int_0^T X_s^2 \,\dd s}
                 \cdot \frac{ \int_0^T X_s\sqrt{Y_s}\,\dd B_s}{\int_0^T X_s^2Y_s\,\dd s}
              + \frac{1}{T}\int_0^T Y_s\,\dd s
               \cdot\frac{ \int_0^T \sqrt{Y_s}\,\dd B_s }{\int_0^T Y_s\,\dd s}}
         {1 - \frac{\left(  \frac{1}{T}\int_0^T X_s \,\dd s \right)^2}{\frac{1}{T}\int_0^T  X_s^2 \,\dd s }} \\
   & =  \frac{ - \EE( X_\infty) \cdot \frac{\EE\left( X_\infty^2 Y_\infty\right)}{\EE( X_\infty^2)} \cdot 0
             + \EE(Y_\infty) \cdot 0 }
            {1 - \frac{(\EE( X_\infty))^2}{\EE(X_\infty^2)} }
    = 0\qquad \text{a.s.}
 \end{align*}
Using that the intersection of two events with probability one is an event with probability one, we have the assertion.
\proofend

\section{Asymptotic behaviour of maximum likelihood estimator}
\label{section_AMLE}

\begin{Thm}\label{Thm_MLE_theta_asymp}
If \ $a>1/2$, \ $b>0$, \ $m\in\RR$, \ $\theta>0$, \ and \ $\PP(Y_0>0)=1$,
 \ then the MLE of \ $\theta$ \ is asymptotically normal, i.e.,
 \begin{align*}
  \sqrt{T}(\widetilde\theta^{\mathrm{MLE}}_T - \theta) \distr \cN\left(0,\frac{1}{ \EE\left(\frac{X_\infty^2}{Y_\infty}\right)}\right)
   \qquad \text{as \ $T\to\infty$,}
 \end{align*}
 where \ $\EE\left(X_\infty^2/Y_\infty\right)$ \ is positive and finite.
\end{Thm}

\noindent{\bf Proof.}
First note that, by \eqref{help52},
 \begin{align}\label{help43}
   \sqrt{T}(\widetilde\theta^{\mathrm{MLE}}_T - \theta)
      = - \frac{\frac{1}{\sqrt{T}}\int_0^T\frac{X_s}{\sqrt{Y_s}}\,\dd B_s}
               {\frac{1}{T}\int_0^T\frac{X_s^2}{Y_s}\,\dd s}
 \qquad \text{a.s.}
 \end{align}
 due to \eqref{help40}.
Recall that, by \eqref{help10},
 \ $\PP\big(\lim_{T\to\infty} \int_0^T \frac{X_s^2}{Y_s}\,\dd s = +\infty\big)=1$.
\ Using Theorem \ref{THM_Zanten} with the following choices
 \begin{align*}
    & d:=1,\qquad
      M_t := \int_0^t \frac{X_s}{\sqrt{Y_s}}\,\dd B_s, \quad t\geq 0,\qquad
      \text{$(\cF_t)_{t\geq 0}$ \ given in Remark \ref{augment}},\\
    & Q(t) :=  \frac{1}{\sqrt{t}}, \quad  t>0,
      \qquad\quad
     \eta:=\sqrt{\EE\left(\frac{X_\infty^2}{Y_\infty}\right)},
 \end{align*}
 we have
 \begin{align*}
   \frac{1}{\sqrt{T}}\int_0^T \frac{X_s}{\sqrt{Y_s}} \,\dd B_s
       \distr  \sqrt{\EE\left(\frac{X_\infty^2}{Y_\infty}\right)}\xi \qquad \text{as \ $T\to\infty$,}
 \end{align*}
 where \ $\xi$ \ is a standard normally distributed random variable.
Then Slutsky's lemma, \eqref{help42} and \eqref{help43} yield the assertion.
\proofend

\begin{Thm}\label{Thm_MLE_theta_m}
If \ $a>1/2$, \ $b>0$, \ $m\in\RR$, \ $\theta>0$, \ and \ $\PP(Y_0>0)=1$,
 \ then the MLE of \ $(a,b,m,\theta)$ \ is asymptotically normal, i.e.,
 \begin{align*}
  \sqrt{T}
    \begin{bmatrix}
      \widehat a^{\mathrm{MLE}}_T - a \\
      \widehat b^{\mathrm{MLE}}_T - b \\
      \widehat m^{\mathrm{MLE}}_T - m \\
      \widehat\theta^{\mathrm{MLE}}_T - \theta
    \end{bmatrix}
    \distr
    \cN_4\left(0,\Sigma^{\mathrm{MLE}}\right) \qquad \text{as \ $T\to\infty$,}
 \end{align*}
 where \ $\cN_4\left(0,\Sigma^{\mathrm{MLE}}\right)$ \ denotes a \ $4$-dimensional normally distribution with
 mean vector \ $0\in\RR^4$ \ and with covariance matrix
 \ $\Sigma^{\mathrm{MLE}} := \diag(\Sigma_1^{\mathrm{MLE}} , \Sigma_2^{\mathrm{MLE}})$ \ with
 blockdiagonals given by
 \begin{gather*}
  \Sigma_1^{\mathrm{MLE}}
  := \frac{1}{\EE\left(\frac{1}{Y_\infty}\right)\EE(Y_\infty) - 1} D_1 , \qquad
  D_1:=\begin{bmatrix}
        \EE(Y_\infty) & 1 \\
        1 & \EE\left(\frac{1}{Y_\infty}\right)
       \end{bmatrix} , \\
  \Sigma_2^{\mathrm{MLE}}
  := \frac{1}{\EE\left(\frac{1}{Y_\infty}\right)
              \EE\left(\frac{X_\infty^2}{Y_\infty}\right)
              - \left(\EE\left(\frac{X_\infty}{Y_\infty}\right)\right)^2} D_2 ,
  \qquad
  D_2:=\begin{bmatrix}
           \EE\left(\frac{X_\infty^2}{Y_\infty}\right)
         & \EE\left(\frac{X_\infty}{Y_\infty}\right) \\
        \EE\left(\frac{X_\infty}{Y_\infty}\right)
         &  \EE\left(\frac{1}{Y_\infty}\right)
       \end{bmatrix} .
 \end{gather*}
\end{Thm}

\noindent{\bf Proof.}
By \eqref{help47a}, \eqref{help47b}, \eqref{help47} and \eqref{help48}, we have
 \begin{align*}
  \sqrt{T}(\widehat{a}_T^{\mathrm{MLE}} - a)
  &=\frac{\frac{1}{T} \int_0^T Y_s \, \dd s \cdot
          \frac{1}{\sqrt{T}} \int_0^T \frac{1}{\sqrt{Y_s}} \, \dd L_s
          - \frac{1}{\sqrt{T}} \int_0^T \sqrt{Y_s} \, \dd L_s}
         {\frac{1}{T} \int_0^T Y_s \, \dd s \cdot
          \frac{1}{T} \int_0^T \frac{1}{Y_s}\,\dd s - 1}
  \qquad \text{a.s.,} \\
  \sqrt{T}(\widehat{b}_T^{\mathrm{MLE}} - b)
  &=\frac{\frac{1}{\sqrt{T}} \int_0^T \frac{1}{\sqrt{Y_s}}\,\dd L_s
          - \frac{1}{T} \int_0^T \frac{1}{Y_s} \, \dd s \cdot
            \frac{1}{\sqrt{T}} \int_0^T \sqrt{Y_s} \, \dd L_s}
         {\frac{1}{T} \int_0^T Y_s \, \dd s \cdot
          \frac{1}{T} \int_0^T \frac{1}{Y_s}\,\dd s - 1}
  \qquad \text{a.s.,} \\
  \sqrt{T}(\widehat{m}_T^{\mathrm{MLE}} - m)
  &=\frac{\frac{1}{T} \int_0^T\frac{X_s^2}{Y_s} \, \dd s \cdot
          \frac{1}{\sqrt{T}} \int_0^T \frac{1}{\sqrt{Y_s}} \, \dd B_s
          - \frac{1}{T} \int_0^T \frac{X_s}{Y_s} \, \dd s \cdot
            \frac{1}{\sqrt{T}} \int_0^T \frac{X_s}{\sqrt{Y_s}} \, \dd B_s}
         {\frac{1}{T} \int_0^T \frac{X_s^2}{Y_s} \, \dd s \cdot
          \frac{1}{T}\int_0^T \frac{1}{Y_s}\,\dd s
          - \left(\frac{1}{T} \int_0^T \frac{X_s}{Y_s} \, \dd s\right)^2}
  \qquad \text{a.s.,} \\
  \sqrt{T}(\widehat\theta_T^{\mathrm{MLE}} - \theta)
  &=\frac{\frac{1}{T}\int_0^T\frac{X_s}{Y_s} \, \dd s \cdot
          \frac{1}{\sqrt{T}} \int_0^T \frac{1}{\sqrt{Y_s}} \, \dd B_s
          - \frac{1}{T} \int_0^T \frac{1}{Y_s} \,\dd s \cdot
            \frac{1}{\sqrt{T}} \int_0^T \frac{X_s}{\sqrt{Y_s}} \, \dd B_s}
         {\frac{1}{T} \int_0^T \frac{X_s^2}{Y_s} \, \dd s
          \cdot\frac{1}{T}\int_0^T \frac{1}{Y_s}\,\dd s
          - \left(\frac{1}{T} \int_0^T \frac{X_s}{Y_s}\,\dd s\right)^2}
  \qquad \text{a.s.,}
 \end{align*}
 due to \eqref{help41a} and \eqref{help41}.
Next, we show that
 \begin{align}\label{help49}
  \frac{1}{\sqrt{T}} M_T
  :=\frac{1}{\sqrt{T}}
    \begin{bmatrix}
     \int_0^T \sqrt{Y_s} \, \dd L_s \\
     \int_0^T \frac{1}{\sqrt{Y_s}} \, \dd L_s \\
     \int_0^T \frac{X_s}{\sqrt{Y_s}} \, \dd B_s \\
     \int_0^T \frac{1}{\sqrt{Y_s}} \, \dd B_s
    \end{bmatrix}
  \distr \eta Z ,
  \qquad \text{as \ $T\to\infty$,}
 \end{align}
 where \ $Z$ \ is a 4-dimensional standard normally distributed random variable and
 \ $\eta$ \ is a non-random $4\times 4$ matrix such that
 \[
   \eta\eta^\top = \diag(D_1, D_2) .
 \]
Here the matrices \ $D_1$ \ and \ $D_2$ \ are positive definite, since their principal minors are positive,
 and \ $\eta$ \ denotes the unique symmetric positive definite square root of \ $\diag(D_1, D_2)$
 (see, e.g., Horn and Johnson \cite[Theorem 7.2.6]{HorJoh}).
\ Indeed, by the absolutely continuity of \ $(Y_\infty, X_\infty)$ \ (see Theorem \ref{Thm_ergodic12}),
 there does not exist constants \ $c,d\in\RR_+$ \ such that \ $\PP(Y_\infty^2=c)=1$ \ and
 \ $\PP(X_\infty^2=d)=1$, \ and, by Cauchy-Schwarz's inequality,
 \begin{align*}
  \EE\left(\frac{1}{Y_\infty}\right)  \EE(Y_\infty) - 1\geq 0
    \qquad \text{and} \qquad
   \EE\left(\frac{X_\infty^2}{Y_\infty}\right) \EE\left(\frac{1}{Y_\infty}\right)
       - \left( \EE\left(\frac{X_\infty}{Y_\infty}\right) \right)^2
      \geq 0,
 \end{align*}
 where equalities would hold if and only if \ $\PP(Y_\infty^2=c)=1$ \ and \ $\PP(X_\infty^2=d)=1$ \ with some constants
  \ $c,d\in\RR_+$, \ respectively.
Note also that the quantity \ $\EE(1/Y_\infty) \EE(Y_\infty) - 1$ \ could have been calculated explicitly, since
 \ $Y_\infty$ \ has Gamma distribution with parameters \ $2a$ \ and \ $2b$.

Let us use Theorem \ref{THM_Zanten} with the choices \ $d=4$, \ $M_t$,
 \ $t\geq 0$, \ defined in \eqref{help49}, \ $(\cF_t)_{t\geq 0}$ \ given in
 Remark \ref{augment}, and
 \ $Q(t) := \diag(t^{-1/2},t^{-1/2},t^{-1/2},t^{-1/2})$, \ $t>0$.
\ We have
 \[
   \langle M \rangle_t
   = \begin{bmatrix}
      \int_0^t Y_s \,\dd s & t & 0 & 0 \\
      t & \int_0^t \frac{1}{Y_s} \, \dd s & 0 & 0 \\
      0 & 0 & \int_0^t \frac{X_s^2}{Y_s} \, \dd s
       & \int_0^t\frac{X_s}{Y_s} \, \dd s \\
      0 & 0 & \int_0^t \frac{X_s}{Y_s} \, \dd s & \int_0^t \frac{1}{Y_s} \,\dd s
     \end{bmatrix} ,
   \qquad t \geq 0 ,
 \]
 where we used the independence of \ $(L_t)_{t\geq 0}$ \ and \ $(B_t)_{t\geq 0}$.
\ Recall that (under the conditions of the theorem) in the proof of Theorem \ref{Thm_MLE_cons_theta_m}
 it was shown that \ $\EE(Y_\infty)<\infty$, \ $\EE(1/Y_\infty)<\infty$,
 \ $\EE(X_\infty/Y_\infty)<\infty$, \ and \ $\EE(X_\infty^2/Y_\infty)<\infty$, \ and, hence
 by Theorem \ref{Thm_ergodic12}, we have
 \[
   Q(t) \langle M \rangle_t Q(t)^\top
   \to \diag(D_1, D_2)
   \qquad \text{as \ $t \to \infty$ \ a.s.}
 \]
Hence, Theorem \ref{THM_Zanten} yields \eqref{help49}.
Then Slutsky's lemma and the continuous mapping theorem yield that
 \begin{align*}
  \sqrt{T}
  \begin{bmatrix}
   \widehat a^{\mathrm{MLE}}_T - a \\
   \widehat b^{\mathrm{MLE}}_T - b \\
   \widehat m^{\mathrm{MLE}}_T - m \\
   \widehat\theta^{\mathrm{MLE}}_T - \theta
  \end{bmatrix}
  \distr
  \diag(A_1, A_2) \eta Z
  \qquad\text{as \ $T \to \infty$,}
 \end{align*}
 where
 \begin{gather*}
  A_1
  := \frac{1}{\EE\left(\frac{1}{Y_\infty}\right)\EE(Y_\infty) - 1} B_1 , \qquad
  B_1:=\begin{bmatrix}
        -1 & \EE(Y_\infty) \\
        -\EE\left(\frac{1}{Y_\infty}\right) & 1
       \end{bmatrix} , \\
  A_2
  := \frac{1}{\EE\left(\frac{1}{Y_\infty}\right)
              \EE\left(\frac{X_\infty^2}{Y_\infty}\right)
              - \left(\EE\left(\frac{X_\infty}{Y_\infty}\right)\right)^2} B_2 ,
  \qquad
  B_2:=\begin{bmatrix}
        -\EE\left(\frac{X_\infty}{Y_\infty}\right)
         & \EE\left(\frac{X_\infty^2}{Y_\infty}\right) \\
        -\EE\left(\frac{1}{Y_\infty}\right)
         & \EE\left(\frac{X_\infty}{Y_\infty}\right)
       \end{bmatrix} .
 \end{gather*}
Using that \ $\eta Z$ \ is a 4-dimensional normally distributed random variable with mean vector zero
 and with covariance matrix \ $\eta\eta^\top = \diag(D_1, D_2)$, \ the covariance matrix of \ $\diag(A_1, A_2) \eta Z$ \ takes the form
 \begin{align*}
  \diag(A_1, A_2) \diag(D_1, D_2) \diag(A_1^\top, A_2^\top)
  = \diag(A_1 D_1 A_1^\top, A_2 D_2 A_2^\top) ,
 \end{align*}
 which yields the assertion.
Indeed,
 \begin{align*}
  B_1 D_1 B_1^\top = \left( \EE\left(\frac{1}{Y_\infty}\right)\EE(Y_\infty) - 1 \right) D_1,
 \end{align*}
 and
 \begin{align*}
  B_2 D_2 B_2^\top = \left( \EE\left(\frac{X_\infty^2}{Y_\infty}\right)\EE\left(\frac{1}{Y_\infty}\right) -
                           \left( \EE\left(\frac{X_\infty}{Y_\infty}\right) \right)^2 \right) D_2.
 \end{align*}
\proofend

\begin{Rem}\label{Rem1}
The asymptotic variance \ $1/\EE\left(X_\infty^2/Y_\infty\right)$ \ of
 \ $\widetilde\theta^{\mathrm{MLE}}_T$ \ in Theorem \ref{Thm_MLE_theta_asymp} is less than
 the asymptotic variance
 \[
     \frac{\EE\left(\frac{1}{Y_\infty}\right)}
          {\EE\left(\frac{1}{Y_\infty}\right)\EE\left(\frac{X_\infty^2}{Y_\infty}\right)
                         - \left(\EE\left(\frac{X_\infty}{Y_\infty}\right)\right)^2}
 \]
 of \ $\widehat\theta^{\mathrm{MLE}}_T$ \ in Theorem \ref{Thm_MLE_theta_m}.
This is in accordance with the fact that \ $\widetilde\theta^{\mathrm{MLE}}_T$ \  is the MLE
 of \ $\theta$ \ provided that the value of the parameter \ $m$ \ is known, which gives extra information,
 so the MLE estimator of \ $\theta$ \ becomes better.
\proofend
\end{Rem}

\section{Asymptotic behaviour of least squares estimator}
\label{section_ALSE}

\begin{Thm}\label{Thm_LSE_theta}
If \ $a>0$, \ $b>0$, \ $m\in\RR$, \ $\theta>0$, \ and \ $\PP(Y_0>0)=1$,
 \ then the LSE of \ $\theta$ \ is asymptotically normal, i.e.,
 \begin{align*}
  \sqrt{T}(\widetilde\theta^{\mathrm{LSE}}_T - \theta) \distr \cN\left(0,\frac{\EE(X_\infty^2Y_\infty)}{ (\EE(X_\infty^2))^2}\right)
   \qquad \text{as \ $T\to\infty$,}
 \end{align*}
 where \ $\EE(X_\infty^2Y_\infty)$ \ and \ $\EE(X_\infty^2)$ \ are given explicitly in Theorem \ref{Thm_ergodic12}.
\end{Thm}

\noindent{\bf Proof.}
First note that, by \eqref{help27},
 \begin{align}\label{help28}
    \sqrt{T}(\widetilde\theta^{\mathrm{LSE}}_T - \theta)
         = - \frac{\frac{1}{\sqrt{T}}\int_0^T X_s\sqrt{Y_s}\,\dd B_s }
                  {\frac{1}{T}\int_0^T X_s^2\,\dd s}
                  \qquad \text{a.s.}
 \end{align}
 due to \eqref{help39}.
Using that \ $\EE(X_\infty^2Y_\infty)$ \ is positive (since \ $X_\infty^2Y_\infty$ \ is non-negative and absolutely
 continuous), by \eqref{help80}, we have
 \begin{align*}
   \PP\left(\lim_{T\to\infty} \int_0^T X_s^2Y_s \,\dd s = +\infty \right)=1.
 \end{align*}
Further, an application of Theorem \ref{THM_Zanten} with the following choices
 \begin{align*}
     &d:=1,\qquad
     M_t := \int_0^t X_s\sqrt{Y_s}\,\dd B_s, \quad t\geq 0,\qquad
     \text{$(\cF_t)_{t\geq 0}$ \ given in Remark \ref{augment},} \\
     & Q(t) := \frac{1}{\sqrt{t}}, \quad t>0,
     \qquad\quad
     \eta:=\sqrt{\EE(X_\infty^2Y_\infty)},
 \end{align*}
 yields that
 \begin{align*}
   \frac{1}{\sqrt{T}}\int_0^T X_s\sqrt{Y_s}\,\dd B_s
       \distr  \sqrt{\EE(X_\infty^2Y_\infty)}\xi \qquad \text{as \ $T\to\infty$,}
 \end{align*}
 where \ $\xi$ \ is a standard normally distributed random variable.
Using again \eqref{help80}, Slutsky's lemma and \eqref{help28}, we get the assertion.
\proofend

\begin{Rem}
The asymptotic variance \ $\EE(X_\infty^2Y_\infty)/(\EE(X_\infty^2))^2$ \ of the LSE
 \ $\widetilde\theta^{\mathrm{LSE}}_T$ \ in Theorem \ref{Thm_LSE_theta} is greater than the asymptotic variance
 \ $1/\EE\left(X_\infty^2/Y_\infty\right)$ \ of the MLE
 \ $\widetilde\theta^{\mathrm{MLE}}_T$ \ in Theorem \ref{Thm_MLE_theta_asymp}, since, by Cauchy and Schwarz's inequality,
 \begin{align*}
   (\EE(X_\infty^2))^2
       = \left(\EE\left(\frac{X_\infty}{\sqrt{Y_\infty}} X_\infty \sqrt{Y_\infty}\right)\right)^2
       < \EE\left(\frac{X_\infty^2}{Y_\infty}\right) \EE\big(X_\infty^2 Y_\infty\big).
 \end{align*}
Note also that using the limit theorem for
 \ $\widetilde\theta^{\mathrm{LSE}}_T$ \ given in Theorem \ref{Thm_LSE_theta},
 one can not give asymptotic confidence intervals for \ $\theta$, \ since the
 variance of the limit normal distribution depends on the unknown parameters
 \ $a$ \ and \ $b$ \ as well.
\proofend
\end{Rem}

\begin{Thm}\label{Thm_LSE_theta_m}
If \ $a>0$, \ $b>0$, \ $m\in\RR$, \ $\theta>0$, \ and \ $\PP(Y_0>0)=1$,
 \ then the LSE of \ $(m,\theta)$ \ is asymptotically normal, i.e.,
 \begin{align*}
  \sqrt{T}
    \begin{bmatrix}
      \widehat m^{\mathrm{LSE}}_T - m \\
      \widehat\theta^{\mathrm{LSE}}_T - \theta \\
    \end{bmatrix}
    \distr
    \cN_2\left(0,\Sigma^{\mathrm{LSE}}\right) \qquad \text{as \ $T\to\infty$,}
 \end{align*}
 where \ $\cN_2\left(0,\Sigma^{\mathrm{LSE}}\right)$ \ denotes a \ $2$-dimensional normally distribution with
 mean vector \ $0\in\RR^2$ \ and with covariance matrix \ $\Sigma^{\mathrm{LSE}}=(\Sigma^{\mathrm{LSE}}_{i,j})_{i,j=1}^2$, \
 where
 \begin{align*}
   &\Sigma^{\mathrm{LSE}}_{1,1}:=\frac{ (\EE(X_\infty))^2\EE(X_\infty^2Y_\infty) - 2\EE(X_\infty)\EE(X_\infty^2)\EE(X_\infty Y_\infty)
                          + (\EE(X_\infty^2))^2\EE(Y_\infty)}
                       {(\EE(X_\infty^2) - (\EE(X_\infty))^2 )^2},\\
   &\Sigma^{\mathrm{LSE}}_{1,2}=\Sigma^{\mathrm{LSE}}_{2,1}:=\frac{\EE(X_\infty)\big( \EE(X_\infty^2Y_\infty) + \EE(X_\infty^2)\EE(Y_\infty) \big)
                                       - \EE(X_\infty Y_\infty)\big( \EE(X_\infty^2) + (\EE(X_\infty))^2 \big)}
                                  {(\EE(X_\infty^2) - (\EE(X_\infty))^2 )^2},\\
   &\Sigma^{\mathrm{LSE}}_{2,2}:= \frac{ \EE(X_\infty^2Y_\infty) - 2\EE(X_\infty)\EE(X_\infty Y_\infty) + (\EE(X_\infty))^2\EE(Y_\infty)}
                        {(\EE(X_\infty^2) - (\EE(X_\infty))^2 )^2}.
 \end{align*}
\end{Thm}

\noindent{\bf Proof.}
By \eqref{help44} and \eqref{help45}, we have
 \begin{align*}
   \sqrt{T}(\widehat m_T^{\mathrm{LSE}} - m)
     =  \frac{-\frac{1}{T}\int_0^T X_s\,\dd s \frac{1}{\sqrt{T}}\int_0^T X_s\sqrt{Y_s}\,\dd B_s
               + \frac{1}{T} \int_0^T X_s^2\,\dd s  \frac{1}{\sqrt{T}}\int_0^T \sqrt{Y_s}\,\dd B_s }
             {\frac{1}{T}\int_0^T X_s^2\,\dd s - \left(\frac{1}{T}\int_0^T X_s\,\dd s\right)^2}
  \qquad \text{a.s.}
 \end{align*}
 and
 \begin{align*}
  \sqrt{T}(\widehat \theta^{\mathrm{LSE}}_T - \theta)
    = \frac{-\frac{1}{\sqrt{T}}\int_0^T X_s\sqrt{Y_s}\,\dd B_s
            + \frac{1}{T}\int_0^T X_s\,\dd s  \frac{1}{\sqrt{T}}\int_0^T \sqrt{Y_s}\,\dd B_s }
             {\frac{1}{T}\int_0^T X_s^2\,\dd s - \left(\frac{1}{T}\int_0^T X_s\,\dd s\right)^2}
  \qquad \text{a.s.}
 \end{align*}
 due to \eqref{help38}.
Next we show that
 \begin{align}\label{help46}
   \left( \frac{1}{\sqrt{T}}\int_0^T X_s\sqrt{Y_s}\,\dd B_s,
          \frac{1}{\sqrt{T}}\int_0^T \sqrt{Y_s}\,\dd B_s \right)
        \distr \Big( (\eta Z)_1, (\eta Z)_2 \Big)
        \quad \text{as \ $T\to\infty$,}
 \end{align}
 where \ $Z$ \ is a 2-dimensional standard normally distributed random variable and
 \ $\eta$ \ is a non-random $2\times 2$ matrix such that
 \[
    \eta\eta^\top = \begin{bmatrix}
                      \EE(X_\infty^2 Y_\infty) & \EE(X_\infty Y_\infty) \\
                      \EE(X_\infty Y_\infty) & \EE(Y_\infty) \\
                    \end{bmatrix}.
 \]
Here the matrix
 \[
  \begin{bmatrix}
     \EE(X_\infty^2 Y_\infty) & \EE(X_\infty Y_\infty) \\
     \EE(X_\infty Y_\infty) & \EE(Y_\infty) \\
  \end{bmatrix}
 \]
 is positive definite, since its principal minors are positive, and \ $\eta$ \ denotes its unique symmetric
 positive definite square root.
Indeed, by the absolutely continuity of \ $(Y_\infty, X_\infty)$ \ (see Theorem \ref{Thm_ergodic12}),
 we have \ $\PP(X_\infty^2 Y_\infty = 0)=0$ \ and, by Cauchy-Schwarz's inequality,
 \ $\EE (X_\infty^2 Y_\infty) \EE(Y_\infty) - \left( \EE(X_\infty Y_\infty) \right)^2 \geq 0$,
 \ where equality would hold if and only if \ $\PP(K X_\infty^2 Y_\infty=L Y_\infty)=1$ \ with some constant \ $K,L\in\RR_+$ \ such that
 \ $K^2+L^2>0$ \ or equivalently (using that \ $\PP(Y_\infty>0)=1$ \ since \ $Y_\infty$ \ has Gamma distribution)
 if and only if \ $\PP(X_\infty^2 = L/K)=1$, \ which leads us to a contradiction
 refereeing to the absolutely continuity of \ $X_\infty$.

Let us use Theorem \ref{THM_Zanten} with the following choices
 \begin{align*}
     &d:=2,\qquad\qquad
     \text{$(\cF_t)_{t\geq 0}$ \ given in Remark \ref{augment},}\\
  & M_t := \begin{bmatrix}
              \int_0^t X_s\sqrt{Y_s}\,\dd B_s \\
              \int_0^t \sqrt{Y_s}\,\dd B_s  \\
            \end{bmatrix},
        \quad t\geq 0,\qquad \quad
  Q(t):=  \begin{bmatrix}
             \frac{1}{\sqrt{t}} & 0 \\
             0 &  \frac{1}{\sqrt{t}} \\
           \end{bmatrix},
            \quad t>0.
      \end{align*}
Then
 \[
    \langle M\rangle_t
       = \begin{bmatrix}
                \int_0^t X_s^2Y_s\,\dd s & \int_0^t X_sY_s\,\dd s \\
                \int_0^t X_sY_s\,\dd s & \int_0^t Y_s\,\dd s \\
          \end{bmatrix},
     \qquad t\geq 0,
 \]
 and hence, by Theorem \ref{Thm_ergodic12} (similarly as detailed in the proof of Theorem \ref{Thm_MLE_theta_m}),
 \[
  Q(t)\langle M\rangle_t  Q(t)^\top
      = \begin{bmatrix}
                \frac{1}{t}\int_0^t X_s^2Y_s\,\dd s & \frac{1}{t}\int_0^t X_sY_s\,\dd s \\
                \frac{1}{t}\int_0^t X_sY_s\,\dd s & \frac{1}{t}\int_0^t Y_s\,\dd s \\
          \end{bmatrix}
       \to \begin{bmatrix}
                      \EE(X_\infty^2 Y_\infty) & \EE(X_\infty Y_\infty) \\
                      \EE(X_\infty Y_\infty) & \EE(Y_\infty) \\
           \end{bmatrix}
       \qquad\text{as \ $t\to\infty$ \ a.s.}
 \]
By \eqref{help79}, \eqref{help80}, Slutsky's lemma and the continuous mapping theorem, we get
 \begin{align*}
  \sqrt{T}
   \begin{bmatrix}
   \widehat m^{\mathrm{LSE}}_T - \theta  \\
    \widehat \theta^{\mathrm{LSE}}_T - \theta \\
   \end{bmatrix}
  \distr \frac{1}{\EE(X_\infty^2) - (\EE(X_\infty))^2}
         \begin{bmatrix}
           -\EE(X_\infty) & \EE(X_\infty^2) \\
            -1 & \EE(X_\infty) \\
         \end{bmatrix}
         \eta Z
   \qquad\text{as \ $T\to\infty$.}
 \end{align*}
Using that \ $\eta Z$ \ is a 2-dimensional normally distributed random variable with mean vector zero
 and with covariance matrix
 \[
    \eta\eta^\top = \begin{bmatrix}
                      \EE(X_\infty^2 Y_\infty) & \EE(X_\infty Y_\infty) \\
                      \EE(X_\infty Y_\infty) & \EE(Y_\infty) \\
                    \end{bmatrix},
 \]
 the covariance matrix of
 \[
   \begin{bmatrix}
     - \EE(X_\infty) & \EE(X_\infty^2) \\
     -1 & \EE(X_\infty ) \\
   \end{bmatrix}
   \eta Z
 \]
 takes the form
 \begin{align*}
   & \begin{bmatrix}
      - \EE(X_\infty) & \EE(X_\infty^2) \\
      -1 & \EE(X_\infty ) \\
     \end{bmatrix}
     \eta \EE(ZZ^\top)\eta^\top
     \begin{bmatrix}
      -\EE(X_\infty) & -1 \\
      \EE(X_\infty^2) & \EE(X_\infty) \\
     \end{bmatrix} \\
   &\quad
     = \begin{bmatrix}
      - \EE(X_\infty) & \EE(X_\infty^2) \\
      -1 & \EE(X_\infty ) \\
     \end{bmatrix}
     \begin{bmatrix}
      \EE(X_\infty^2 Y_\infty) & \EE(X_\infty Y_\infty) \\
      \EE(X_\infty Y_\infty) & \EE(Y_\infty) \\
     \end{bmatrix}
     \begin{bmatrix}
      -\EE(X_\infty) & -1 \\
      \EE(X_\infty^2) & \EE(X_\infty) \\
     \end{bmatrix} \\
   &\quad
     = \begin{bmatrix}
         \EE(X_\infty^2)\EE(X_\infty Y_\infty) - \EE(X_\infty) \EE(X_\infty^2 Y_\infty)
                  & \EE(X_\infty^2)\EE(Y_\infty) - \EE(X_\infty)\EE(X_\infty Y_\infty)\\
         \EE(X_\infty)\EE(X_\infty Y_\infty) - \EE(X_\infty^2 Y_\infty)
                  & \EE(X_\infty)\EE(Y_\infty) - \EE(X_\infty Y_\infty) \\
        \end{bmatrix}\\
   &\phantom{\quad=\;}
      \times \begin{bmatrix}
      -\EE(X_\infty) & -1 \\
      \EE(X_\infty^2) & \EE(X_\infty) \\
     \end{bmatrix},
 \end{align*}
 which yields the assertion.
\proofend

\begin{Rem}
Using the explicit forms of the mixed moments given in (iii) of Theorem \ref{Thm_ergodic12},
 one can check that the asymptotic variance \ $\EE(X_\infty^2Y_\infty)/(\EE(X_\infty^2))^2$ \ of
 \ $\widetilde\theta^{\mathrm{LSE}}_T$ \ in Theorem \ref{Thm_LSE_theta} is less than the asymptotic variance
 \ $\Sigma^{\mathrm{LSE}}_{1,1}$ \ of \ $\widehat\theta^{\mathrm{LSE}}_T$ \ in Theorem \ref{Thm_LSE_theta_m}.
This can be interpreted similarly as in Remark \ref{Rem1}.
Note also that using the limit theorem for
 \ $(\widehat m^{\mathrm{LSE}}_T,\widehat \theta^{\mathrm{LSE}}_T)$ \ given
 in Theorem \ref{Thm_LSE_theta_m}, one can not give asymptotic confidence
 regions for \ $(m, \theta)$, \ since the variance matrix of the limit normal
 distribution depends on the unknown parameters
 \ $a$ \ and \ $b$ \ as well.
\proofend
\end{Rem}

\vspace*{3mm}

\appendix

\vspace*{5mm}

\noindent{\bf\Large Appendix}

\section{Radon-Nykodim derivatives for certain diffusions}

We consider the SDEs
 \begin{align}\label{SDE_1}
  \dd \xi_t &= (A \xi_t + a) \, \dd t + \sigma(\xi_t) \,\dd W_t,
  \qquad t\in\RR_+,\\
  \dd \eta_t &= (B \eta_t + b) \, \dd t + \sigma(\eta_t) \,\dd W_t ,
  \qquad t\in\RR_+, \label{SDE_2}
 \end{align}
 with the same initial values \ $\xi_0 = \eta_0$, \ where
 \ $A, B \in \RR^{2\times2}$, \ $a, b \in \RR^2$,
 \ $\sigma : \RR^2 \to \RR^{2\times2}$ \ is a Borel measurable function, and
 \ $(W_t)_{t\in\RR_+}$ \ is a two-dimensional standard Wiener process.
Suppose that the SDEs \eqref{SDE_1} and \eqref{SDE_2} admit pathwise unique
 strong solutions.
Let \ $\PP_{(A,a)}$ \ and \ $\PP_{(B,b)}$ \ denote the probability measures on
 the measurable space \ $(\cC(\RR_+, \RR^2),\cB(\cC(\RR_+, \RR^2)))$ \ induced
 by the processes \ $(\xi_t)_{t\in\RR_+}$ \ and \ $(\eta_t)_{t\in\RR_+}$,
 \ respectively.
Here \ $\cC(\RR_+, \RR^2)$ \ denotes the set of continuous \ $\RR^2$-valued
 functions defined on \ $\RR_+$, \ $\cB(\cC(\RR_+, \RR^2))$ \ is the Borel
 \ $\sigma$-algebra on it, and we suppose that the space
 \ $(\cC(\RR_+, \RR^2),\cB(\cC(\RR_+, \RR^2)))$ \ is endowed with the natural
 filtration \ $(\cA_t)_{t\in\RR_+}$, \ given by
 \ $\cA_t:=\varphi_t^{-1}(\cB(\cC(\RR_+, \RR^2)))$, \ where
 \ $\varphi_t:\cC(\RR_+, \RR^2)\to \cC(\RR_+, \RR^2)$ \ is the mapping
 \ $\varphi_t(f)(s):=f(t\wedge s)$, \ $s \in \RR_+$.
\ For all \ $T>0$, \ let \ $\PP_{(A,a),T}:=\PP_{(A,a)}\vert_{\cA_T}$ \ and
 \ $\PP_{(B,b),T}:=\PP_{(B,b)}\vert_{\cA_T}$ \ be the restrictions of
 \ $\PP_{(A,a)}$ \ and \ $\PP_{(B,b)}$ \ to \ $\cA_T$, \ respectively.

From the very general result in Section 7.6.4 of Liptser and Shiryaev
 \cite{LipShiI}, one can deduce the following lemma.

\begin{Lem}\label{LEMMA_LipShi}
Let \ $A, B \in \RR^{2\times2}$, \ $a, b \in \RR^2$, \ and let
 \ $\sigma : \RR^2 \to \RR^{2\times2}$ \ be a Borel measurable function.
Suppose that the SDEs \eqref{SDE_1} and \eqref{SDE_2} admit pathwise unique
 strong solutions.
Let \ $\PP_{(A,a)}$ \ and \ $\PP_{(B,b)}$ \ denote the probability measures
 induced by the unique strong solutions of the SDEs \eqref{SDE_1} and
 \eqref{SDE_2} with the same initial value \ $\xi_0 = \eta_0$, \ respectively.
Suppose that \ $\PP(\exists\,\sigma(\xi_t)^{-1})= 1$ \ and
 \ $\PP(\exists\,\sigma(\eta_t)^{-1})= 1$ \ for all \ $t \in \RR_+$, \ and
 \[
   \PP\left( \int_0^t \left[ (A \xi_s + a)^\top
                             (\sigma(\xi_s) \sigma(\xi_s)^\top)^{-1}
                             (A \xi_s + a)^\top
                             + (B \xi_s + b)^\top
                               (\sigma(\xi_s) \sigma(\xi_s)^\top)^{-1}
                               (B \xi_s + b)^\top \right] \dd s < \infty \right)
   = 1
 \]
 and
 \[
   \PP\left( \int_0^t \left[ (A \eta_s + a)^\top
                             (\sigma(\eta_s) \sigma(\eta_s)^\top)^{-1}
                             (A \eta_s + a)^\top
                             + (B \eta_s + b)^\top
                               (\sigma(\eta_s) \sigma(\eta_s)^\top)^{-1}
                               (B \eta_s + b)^\top \right] \dd s < \infty \right)
   = 1
 \]
 for all \ $t \in \RR_+$.
\ Then for all \ $T > 0$, \ the probability measures \ $\PP_{(A,a),T}$ \ and
 \ $\PP_{(B,b),T}$ \ are absolutely continuous with respect to each other, and
 the Radon-Nykodim derivative of \ $\PP_{(A,a),T}$ \ with respect to
 \ $\PP_{(B,b),T}$ \ (so called likelihood ratio) takes the form
 \begin{align*}
  L^{(A,a),(B,b)}_T((\xi_s)_{s\in[0,T]})
  =\exp\Bigg\{&\int_0^T
                (A \xi_s + a - B \xi_s - b)^\top
                (\sigma(\xi_s) \sigma(\xi_s)^\top)^{-1}
                \, \dd \xi_s \\
              &-\frac{1}{2}
                \int_0^T
                 (A \xi_s + a - B \xi_s - b)^\top
                (\sigma(\xi_s) \sigma(\xi_s)^\top)^{-1}
                 (A \xi_s + a + B \xi_s + b)^\top \dd s
       \Bigg\} .
 \end{align*}
\end{Lem}

We call the attention that conditions (4.110) and (4.111) are also required
 for Section 7.6.4 in Liptser and Shiryaev \cite{LipShiI}, but the Lipschitz
 condition (4.110) in Liptser and Shiryaev \cite{LipShiI} does not hold in
 general for the SDEs \eqref{SDE_1} and \eqref{SDE_2}.
However, we can use formula (7.139) in Liptser and Shiryaev \cite{LipShiI},
 since they use their conditions (4.110) and (4.111) only in order to ensure
 the SDE they consider in Section 7.6.4 has a pathwise unique strong solution
 (see, the proof of Theorem 7.19 in Liptser and Shiryaev \cite{LipShiI}),
 which is supposed in Theorem \ref{LEMMA_LipShi}.

\section*{Acknowledgements}

We are grateful for both referees for their several valuable comments yielding
 an improvement of the manuscript.

\indent
   {\sc M\'aty\'as Barczy},
   Faculty of Informatics, University of Debrecen,
   Pf.~12, H--4010 Debrecen, Hungary.
   Tel.: +36-52-512900, Fax: +36-52-512996,
   e--mail: barczy.matyas@inf.unideb.hu

\indent
    {\sc Leif D\"oring},
    Institut f\"ur Mathematik,
    Universit\"at Z\"urich,
    Winterthurerstrasse 190,
    CH-8057 Z\"urich,
    Switzerland.
    Tel.: +41-(0)44-63 55892, Fax: +41-(0)44-63 55705,\\
    e--mail: leif.doering@googlemail.com

\indent
     {\sc Zenghu Li},
     School of Mathematical Sciences, Beijing Normal University,
     Beijing 100875,  People's Republic of China.
     Tel.: +86-10-58802900, Fax: +86-10-58808202,
     e--mail: lizh@bnu.edu.cn

\indent
      {\sc Gyula Pap},
      Bolyai Institute, University of Szeged,
      Aradi v\'ertan\'uk tere 1, H--6720 Szeged, Hungary.
      Tel.: +36-62-544033, Fax: +36-62-544548,
      e--mail: papgy@math.u-szeged.hu

\end{document}